\documentclass[a4paper,10pt,smallextended]{amsart}
\usepackage{verbatim,latexsym,exscale,amssymb,amsthm,amsmath,amsfonts,mathrsfs}
\usepackage[all]{xy}

\numberwithin{equation}{section}

\xyoption{ps}
\SelectTips{cm}{}

\newdir{ >}{{}*!/-6pt/@{>}}

\newdir{ >}{{}*!/-5pt/\dir{>}}

\newcommand{\id}{\mathrm{id}}

\newcommand{\loc}{\mathrm{loc}}

\newcommand{\tcat}{\mathcal{T}}

\newcommand{\spec}{\mathop{\mathrm{Spec}}}
\newcommand{\specgr}{\mathop{\mathrm{Spec^h}}}

\newcommand{\spc}{\mathop{\mathrm{Spc}}}
\newcommand{\Spc}{\mathop{\mathrm{Spc}}}

\newcommand{\supp}{\mathrm{supp}}
\newcommand{\Supp}{\mathrm{Supp}}
\newcommand{\Supptot}{\mathrm{Supp_{tot}}}

\newcommand{\local}{\mathrm{Fib}}

\newcommand{\Res}[2]{\mathop{\mathrm{Res}_{#1}^{#2}}}
\newcommand{\Ind}[2]{\mathop{\mathrm{Ind}_{#1}^{#2}}}

\newcommand{\Hom}{\mathrm{Hom}}
\newcommand{\inthom}{\mathrm{\underline{Hom}}}

\newcommand{\End}{\mathrm{End}}

\newcommand{\Ann}{\mathrm{Ann}}

\newcommand{\op}{\mathrm{op}}

\newcommand{\cone}{\mathop{\mathrm{cone}}}

\newcommand{\unit}{\mathbf{1}}
\newcommand{\obj}{\mathop{\mathrm{obj}}}

\newcommand{\Ker}{\mathsf{Ker}}

\newcommand{\Image}{\mathsf{Im}}

\newcommand{\Cstar}{\mathrm{C}^*\!}

\newcommand{\Ext}{\mathrm{Ext}}
\newcommand{\Tor}{\mathrm{Tor}}

\newcommand{\KK}{\mathsf{KK}}
\newcommand{\CI}{\mathsf{CI}}
\newcommand{\CC}{\mathsf{CC}}

\newcommand{\Boot}{\mathsf{Boot}}

\newcommand{\proj}{\mathop{\mathrm{Proj}}}

\newcommand{\Modules}{\textrm{-}\mathsf{Mod}}
\newcommand{\VS}{\textrm{-}\mathsf{VS}}

\newcommand{\stmod}{\mathsf{stmod}}
\newcommand{\StMod}{\mathsf{StMod}}

\newcommand{\Projobj}{\mathsf{Proj}}

\newcommand{\D}{\mathsf{D}}

\newcommand{\Ab}{\mathsf{Ab}}

\newcommand{\Osheaf}{\mathcal{O}}

\newcommand{\Z}{\mathbb{Z}}
\newcommand{\R}{\mathbb{R}}
\newcommand{\C}{\mathbb{C}}
\newcommand{\Q}{\mathbb{Q}}

\theoremstyle{definition}

\newtheorem{defi}[equation]{Definition}

\newtheorem{convention}[equation]{Convention}
\newtheorem*{conv*}{Conventions}
\newtheorem*{ack}{Acknowledgements}

\theoremstyle{theorem}

\newtheorem{thm}[equation]{Theorem}
\newtheorem{lemma}[equation]{Lemma}
\newtheorem{thm-defi}[equation]{Theorem-Definition}
\newtheorem{prop}[equation]{Proposition}
\newtheorem{cor}[equation]{Corollary}
\newtheorem{conj}[equation]{Conjecture}
\newtheorem*{thm*}{Theorem}
\newtheorem*{lemma*}{Lemma}
\newtheorem*{cor*}{Corollary}
\newtheorem*{conj*}{Conjecture}

\theoremstyle{remark}

\newtheorem{notation}[equation]{Notation}
\newtheorem{terminology}[equation]{Terminology}
\newtheorem{remark}[equation]{Remark}
\newtheorem*{remark*}{Remark}
\newtheorem{remarks}[equation]{Remarks}
\newtheorem{construction}[equation]{Construction}
\newtheorem{example}[equation]{Example}

\def\hyph{-\penalty0\hskip0pt\relax}

\begin{document}
\title{Tensor triangular geometry and  $KK$-theory  
 }
\author{Ivo Dell'Ambrogio} \date{}
\maketitle

\begin{abstract}
This is a first foray of \emph{tensor triangular geometry} \cite{balmer_prime} into the realm of bivariant topological $K$-theory. As a motivation, we first establish a connection between the Balmer spectrum $\spc(\KK^G)$ and a strong form of the Baum-Connes conjecture with coefficients for the group~$G$, as studied in \cite{meyernest-bc}. We then turn to more tractable categories, namely, the thick triangulated subcategory $\mathcal K^G\subset \KK^G$ and the localizing subcategory $\tcat^G\subset \KK^G$ generated by the tensor unit~$\C$. 
For $G$ finite, we construct for the objects of $\tcat^G$ a support theory in $\spec(R(G))$ with good properties. We see as a consequence that $\spc(\mathcal K^G)$ contains a copy of the Zariski spectrum $\spec(R(G))$ as a retract, where $R(G)=\End_{\KK^G}(\mathbb C)$ is the complex character ring of~$G$. Not surprisingly, we find that $\spc(\mathcal K^{\{1\}})\simeq \spec(\Z)$.
\end{abstract}

\tableofcontents

\section{Introduction}

Let $G$ be a second countable locally compact Hausdorff group, and let $\KK^G$ denote the $G$-equivariant Kasparov category of separable $G$-$\Cstar$-algebras (\cite{kasparov} \cite{mayer-cat}). As shown in \cite{meyernest-bc}, $\KK^G$ is naturally equipped with the structure of a tensor triangulated category (Def.\ \ref{defi:tensor_tr}).  
This means that we are in the domain of \emph{tensor triangular geometry}.
In particular, the (essentially small) category $\KK^G$ has a spectrum $\spc(\KK^G)$, as defined by Paul Balmer \cite{balmer_prime} (see Def.\ \ref{defi:spectrum} below). If $H\leq G$ is a subgroup, the restriction functor $\Res{G}{H}:\KK^G\to\KK^H$ induces a continuous map $(\Res{G}{H})^*:\spc(\KK^H)\to \spc(\KK^G)$. Then

\begin{thm} \label{thm:covering_intro}
Assume that $G$ is such that $\spc(\KK^G)=\bigcup_H\big(\Res{G}{H}\big)^*\big(\spc(\KK^H)\big)$, where $H$ runs through all compact subgroups of~$G$. Then $G$ satisfies the Baum-Connes conjecture for every functor on $\KK^G$ and any coefficient algebra~$A\in\KK^G$. 
\end{thm}

This is proved in $\S$\ref{sec:BC}, where the reader may also find the precise meaning of the conclusion. Now, we do not know yet if the above fact may provide a way of proving Baum-Connes. For one thing, we still don't know of a single non-compact group satisfying the above covering hypothesis. But the result looks intriguing, and it suggests that further \emph{geometric} inquiry in this context will be fruitful. 

As a first step in this direction, we turn to the subcategories $\tcat^G:=\langle \unit\rangle_{\loc}\subset \KK^G$ and $\mathcal K^G:= \langle \unit\rangle \subset \KK^G$, that is, the localizing, respectively the thick triangulated subcategory generated by the tensor unit $\unit= \C\in \KK^G$.  Moreover, we restrict our attention to the much better understood case when the group $G$ is compact or even finite. 
Then the endomorphism ring $\End(\unit)$ of the $\otimes$-unit can be identified with the complex representation ring $R(G)$ of the compact group, which is known to be noetherian if $G$ is a Lie group (e.g.\ finite); see \cite{segal-spec}.
Note that $\mathcal K^G=(\tcat^G)_c$ is the subcategory of compact objects in~$\tcat^G$ (see $\S$\ref{subsec:Br_Bou} and $\S$\ref{subsec:cats}).
When $G=\{1\}$ is trivial, $\Boot:=\tcat^G$ is better known as the ``Bootstrap'' category of separable $\Cstar$-algebras. We will prove in~$\S$\ref{subsec:boot}:

\begin{thm}\label{thm:boot_intro}
There is a canonical homeomorphism $\spc(\Boot_c)\simeq \spec(\Z)$.
\end{thm}

The latter statement generalizes naturally as follows:

\begin{conj} \label{conj:spec_intro}
For every finite group $G$, the natural map $\rho_{\mathcal K^G}:\spc(\mathcal K^G)\to\spec(R(G))$ (see \cite{balmer_spec2} or $\S$\ref{subsec:split} below) is a homeomorphism.
\end{conj}

If true, this would show that, in yet another branch of mathematics, an object of classical interest (here: the spectrum of the complex representation ring of a finite group) can be recovered as the Balmer spectrum of a naturally arising $\otimes$-triangulated category. 
We have some interesting facts that suggest a positive answer. Namely:

\begin{thm}[{Thm.\ \ref{thm:sigma} and Prop.\ \ref{prop:split}}] \label{thm:gen_supp_intro}
Let $G$ be a finite group. Then there exists an assignment $\sigma_G:\obj(\tcat^G)\to 2^{\spec(R(G))}$ from objects of $\tcat^G$ to subsets of the spectrum enjoying the following properties:
\begin{itemize}
\item[(a)] $\sigma_G(0)=\emptyset$ and $\sigma_G(\unit)=\spec(R(G))$.
\item[(b)] $\sigma_G(A\oplus B)= \sigma_G(A)\cup \sigma_G(B)$.
\item[(c)] $\sigma_G(TA)=\sigma_G(A)$.
\item[(d)] $\sigma_G(B)\subset\sigma_G(A)\cup \sigma_G(C)$ for every exact triangle $A\to B\to C\to TA$.
\item[(e)] $\sigma_G(A\otimes B)=\sigma_G(A)\cap\sigma_G(B)$.
\item[(f)] $\sigma_G(\coprod_iA_i)=\bigcup_i\sigma_G(A_i)$.
\item[(g)] if $A\in\mathcal K^G$, then $\sigma_G(A)$ is a closed subset of $\spec(R(G))$.
\end{itemize}
Here $A,B\in \tcat^G$ are any objects and $\coprod_iA_i$ any coproduct in~$\tcat^G$.
In particular, the restriction of $\sigma_G$ to $\mathcal K^G$ is a support datum in the sense of Balmer~\cite{balmer_prime} (see $\S$\ref{subsec:spec} below), so it induces a canonical map $f_G:\spec(R(G))\to \spc(\mathcal K^G)$. This map is topologically split injective; indeed, it provides a continuous section of~$\rho_{\mathcal K^G}$.
\end{thm}

\begin{remark*}
In the course of proving Theorem \ref{thm:gen_supp_intro} we construct, for $G$ compact, a well-behaved `localization of $\tcat^G$ at a prime $\mathfrak p\in \spec(R(G))$', written $\tcat^G_{\mathfrak p}\subset \tcat^G$ (see $\S$\ref{sec:locKth}). It follows for instance that there is a functor $L_{\mathfrak p}:\KK^G\to \tcat^G_{\mathfrak p}$ together with a natural isomorphism  $K^G_*(L_{\mathfrak p}A)\simeq K^G_*(A)_{\mathfrak p}$, for all $A\in \KK^G$ (Cor.\ \ref{cor:local_approx}). 
\end{remark*}

We believe Theorem \ref{thm:gen_supp_intro} provides evidence for Conjecture \ref{conj:spec_intro} because of the following more general result in tensor triangular geometry, which is of independent interest (see Theorem \ref{thm:abstract} below).

\begin{thm} \label{thm:abstract_intro}
Let $\tcat$ be a compactly generated $\otimes$-tri\-angulat\-ed category\footnote{See Convention \ref{conv:cptly_gen} below for the precise (modest) hypotheses we are making here. We require in particular that compact objects form a \emph{tensor} triangulated subcategory~$\tcat_c$.}. Let $X$ be a spectral topological space (such as the Zariski spectrum of a commutative ring -- see Remark \ref{rem:spectral}), and let $\sigma: \obj(\tcat)\to 2^X$ be a function assigning to every object of $\tcat$ a subset of~$X$. Assume that the pair $(X,\sigma)$ satisfies the following ten axioms:
\begin{itemize}
\item[(S0)] $\sigma(0)=\emptyset$.
\item[(S1)] $\sigma(\unit)=X$.
\item[(S2)] $\sigma(A\oplus B)= \sigma(A)\cup \sigma(B)$ (really, this is redundant because of (S6) below).
\item[(S3)] $\sigma(TA)=\sigma(A)$.
\item[(S4)] $\sigma(B)\subset\sigma(A)\cup \sigma(C)$ for every distinguished triangle  $A\to B\to C\to TA$.
\item[(S5)] $\sigma(A\otimes B)=\sigma(A)\cap\sigma(B)$ for every compact $A\in \tcat_c$ and arbitrary $B\in \tcat$. 
\item[(S6)] $\sigma(\coprod_i A_i)= \bigcup_i \sigma(A_i)$ for every small family $\{A_i\}_i\subset \tcat$ of objects.
\item[(S7)] $\sigma(A)$ is closed in $X$ with quasi-compact complement $X\smallsetminus \sigma(A)$ for all  $A\in \tcat_c$.
\item[(S8)] For every closed subset $Z\subset X$ with quasi-compact open complement, there exists a compact object $A\in \tcat_c$ with $\sigma(A)=Z$.
\item[(S9)] $\sigma(A)=\emptyset $ implies $ A\simeq 0$.
\end{itemize}
Then the restriction of $(X,\sigma)$ to $\tcat_c$ is a classifying support datum; in particular, the induced canonical map $X\to \spc(\tcat_c)$ is a homeomorphism (see Thm.\ \ref{thm:abstract_class}).
\end{thm}

\begin{remark} We note that the latter theorem has also been announced by Julia Pevtsova and Paul Smith. 
It specializes to the classification of thick tensor ideals in the stable category $\stmod(kG)$ of modular representation theory, due to  Benson, Carlson and Rickard \cite{bcr}  (see Example \ref{ex:class_bcr} below). Indeed, our proof is an abstract version of their \cite[Theorem~3.4]{bcr}. 
\end{remark}

As concerns us here, our hope is to apply Theorem \ref{thm:abstract_intro} to the category $\tcat:=\tcat^G$ (so that $\tcat_c=\mathcal K^G$) for a finite group~$G$, choosing $\sigma$ to be the assignment $\sigma_G$ in Theorem\ \ref{thm:gen_supp_intro}; note that it follows from the first part of the theorem that $\sigma_G$ satisfies conditions (S0)-(S7). At least for $G=\{1\}$, axioms (S8) and (S9) are also satisfied and therefore we obtain Theorem \ref{thm:boot_intro} from Theorem \ref{thm:abstract_intro}. We don't know yet if the same strategy also works in general, i.e., we don't know if (S8) and (S9) also hold when $G$ is non-trivial (we have some clues that this might be the case, but they are too sparse to be mentioned here).

More abstractly, in $\S$\ref{subsec:enter} we examine condition (S8) (and also (S7)) in relation to the endomorphism ring of the tensor unit~$\unit$. As a payoff, we then show in $\S$\ref{subsec:comparison} how to use Theorem \ref{thm:abstract_intro} in order to compare Balmer's universal support with that of Benson, Iyengar and Krause \cite{bik} in the situation where both are defined.

In a sequel to this article, we intend to study the spectrum of ``finite noncommutative $G$-CW-complexes'' for a finite group~$G$, that is, of the triangulated subcategory of $\KK^G$ generated by all $G$-$\Cstar$-algebras $C(G/H)$ with $H\leq G$ a subgroup.

\begin{conv*}
 If $F:\mathcal A\to \mathcal B$ is an additive functor, we denote by $\mathsf{Im}(F)\subset \mathcal B$ the essential image of~$F$, i.e., the full subcategory of $\mathcal B$ of those objects isomorphic to $F(A)$ for some $A\in \mathcal A$; by $\mathsf{Ker}(F):=\{{A\in \mathcal A} \mid {F(A)\simeq 0}\}$ we denote its kernel on objects, and by $\mathsf{ker}(F):=\{f\in \mathrm{Mor}(\mathcal A)\mid F(f)=0\}$ its kernel on morphisms.
The translation functor in all triangulated categories is denoted by~$T$. Triangulated subcategories are always full and closed under isomorphic objects.
\end{conv*}

\begin{ack}
This work was done during my PhD thesis under the supervision of Paul Balmer. I wish to thank him for his interest and generosity. I am very grateful to Amnon Neeman for spotting two mistakes in a previous version of this paper.
\end{ack}

\section{Triangular preliminaries}

\subsection{Brown representability and Bousfield localization} 
\label{subsec:Br_Bou}

The material of this section, originated in stable homotopy and generalized to triangulated categories by Amnon Neeman in a series of papers, is now standard. However we shall have to use a slight variation of the definitions and results.  
Namely, we fix an uncountable regular cardinal number~$\alpha$, and consider variants of the usual notions that are relative to this cardinal. (Later on, in our applications we shall only need the case $\alpha=\aleph_1$.) We use subscripts as in ``dummyword${}_{\alpha}$'', because the prefixed notation ``$\alpha$\hyph dummyword''  has already found different uses. Throughout, $\tcat$ will be a triangulated category admitting arbitrary \emph{small${}_{\alpha}$ coproducts}, i.e., coproducts indexed by sets $I$ of cardinality $|I|<\alpha$. In general, we shall say that a set $S$ is \emph{small${}_{\alpha}$} if $|S|<\alpha$.

\begin{defi} \label{defi:cpt}
An object $A$ of $\tcat$ is \emph{compact${}_{\alpha}$} if $\Hom_{\tcat}(A,\,?)$ commutes with small${}_{\alpha}$ coproducts, and if moreover $|\Hom_{\tcat}(A,B)|<\alpha$ for every~$B\in \tcat$. We write $\tcat_c$ for the full subcategory of compact${}_{\alpha}$ objects of~$\tcat$. A set of objects $\mathcal G\subset \tcat$ \emph{generates}~$\tcat$ if for all $A\in \tcat$ the following implication holds:
\begin{eqnarray*}
\Hom_{\tcat}(G,A)\simeq 0 \textrm{ for all }G\in \mathcal G\;\; \Rightarrow \;\; A\simeq 0.
\end{eqnarray*}
We say that $\tcat$ is \emph{compactly${}_{\alpha}$ generated} if there is a small${}_{\alpha}$ set $\mathcal G\subset \tcat$ of compact${}_{\alpha}$ objects generating the category.
If $\mathcal E\subset \tcat$ is some class of objects, we write $\langle \mathcal E\rangle_{\loc}$ for the smallest localizing${}_{\alpha}$ subcategory of $\tcat$ containing $\mathcal E$, where \emph{localizing${}_{\alpha}$} means triangulated and closed under the formation of small${}_{\alpha}$ coproducts in~$\tcat$. We will reserve the notation $\langle \mathcal E\rangle $ for the thick triangulated subcategory of $\tcat$ generated by~$\mathcal E$. Note that $\langle \mathcal E\rangle_{\loc}$ is automatically thick, as is every triangulated category with arbitrary countable coproducts, by a well-known argument.

\end{defi}

It was first noticed in \cite{meyernest-bc} that these definitions\footnote{beware that our terminology is slightly changed from that in \emph{loc.\ cit.}} allow the following $\alpha$-relative version of Neeman's Brown representability for cohomological functors, simply by verifying that the usual proof (\cite[Thm.\ 3.1]{neemanBou}) only needs the formation of small${}_{\alpha}$ coproducts in $\tcat$ and never requires bigger ones. 

\begin{thm}[Brown representability] \label{thm:brown} 
Let $\tcat$ be compactly${}_{\alpha}$ generated, with $\mathcal G$ a generating set.  Then a functor $F:\tcat^{\op} \to \Ab$ is representable if and only if it is homological, it sends small${}_{\alpha}$ coproducts in $\tcat$ to products of abelian groups and if moreover $|F(A)|<\alpha$ for all $A\in \mathcal G$ (or equivalently, for all compact${}_{\alpha}$ objects $A\in \tcat_c$). \qed
\end{thm}

As in the case of a \emph{genuine} compactly generated category (i.e., when $\alpha=$ cardinality of a proper class), one obtains from the techniques of the proof the following characterization:

\begin{cor} \label{cor:brown}
For a triangulated category $\tcat$ with arbitrary small${}_{\alpha}$ coproducts, the following are equivalent:
\begin{itemize}
\item[(i)] $\tcat$ is compactly${}_{\alpha}$ generated.
\item[(ii)] $\tcat = \langle \mathcal G\rangle_{\loc}$ for some small${}_{\alpha}$ subset $\mathcal G\subset \tcat_c$ of compact${}_{\alpha}$ objects.
\item[(iii)] $\tcat = \langle \tcat_c\rangle_{\loc}$ and $\tcat_c$ is essentially small${}_{\alpha}$ (by which of course we mean that $\tcat_c$ has a small${}_{\alpha}$ set of isomorphism classes of objects).
\end{itemize}
\end{cor}

\begin{cor} \label{cor:cpt_thick}
Thus, for every small${}_{\alpha}$ subset $\mathcal S\subset \tcat_c$ there is a compactly${}_{\alpha}$ generated localizing${}_{\alpha}$ subcategory $\mathcal L=\langle \mathcal S\rangle_{\loc}\subset \tcat$. Its compact${}_{\alpha}$ objects are given by $\mathcal L_c= \tcat_c\cap \mathcal L= \langle \mathcal S\rangle$.\qed
\end{cor}

\begin{notation}
Let $\mathcal E$ be a class of objects in $\tcat$ closed under translations. We write
\begin{eqnarray*}
\mathcal E^{\perp}&:=&\{A\in \tcat\mid \Hom(E,A)\simeq 0 \textrm{ for all } E\in \mathcal E\} \\
{}^{\perp}\mathcal E&:=&\{A\in \tcat\mid \Hom(A,E)\simeq 0 \textrm{ for all } E\in \mathcal E\} 
\end{eqnarray*}
For two collections $\mathcal E,\mathcal F\subset \tcat$ of objects we write $\mathcal E\perp \mathcal F$ to mean that $\Hom(E,F)\simeq 0$ for all $E\in \mathcal E$ and $F\in\mathcal F$.
\end{notation}

The following proposition collects well-known facts related to Bousfield localization, which we recall in order to fix notation (see e.g.\ \cite[{$\S$9}]{neemanTr}, \cite[{$\S$2.6}]{meyernest-bc}).

\begin{prop}[Bousfield localization] \label{prop:compl}
Let $\tcat$ be a triangulated category, and let $\mathcal L,\mathcal R\subset \tcat$ be thick subcategories satisfying the following condition:
\begin{itemize}
\item[($*$)]  $\mathcal L \perp \mathcal R$ and for every $A\in \tcat$ there exists a distinguished triangle $A'\to A\to A'' \to TA'$ with $A'\in \mathcal L$ and $A''\in \mathcal R$.
\end{itemize}
Then the triangle in ($*$) is unique up to unique isomorphism and is functorial in~$A$. Moreover, the resulting functors $L: A\mapsto A'$ and $R: A\mapsto A''$ and morphisms $\lambda: L\to \id_{\tcat}$ and $\rho: \id_{\tcat}\to R$ enjoy the following properties:
\begin{itemize}
\item[(a)] $\lambda_A:LA\to A$ is the terminal morphism to $A$ from an object of~$\mathcal L$. Dually, $\rho_A:A\to RA$ is initial among morphisms from $A$ to an object of~$\mathcal R$.
\item[(b)] $\mathcal R=\mathcal L^{\perp}$ and $\mathcal L={}^{\perp}\mathcal R$. In particular, $\mathcal L$ and $\mathcal R$ determine each other.
\item[(c)] $\mathcal L$ is a coreflective subcategory of $\tcat$. Dually, $\mathcal L^{\perp}$ is a reflective subcategory.
\item[(d)] The composition $\mathcal L\hookrightarrow \tcat \to \tcat/\mathcal L^{\perp}$ is an equivalence identifying the right adjoint of the inclusion $\mathcal L\hookrightarrow \tcat$ with the Verdier quotient $\tcat\to \tcat/\mathcal L^{\perp}$. Dually, the composition $\mathcal L^{\perp}\hookrightarrow \tcat \to \tcat/\mathcal L$ is an equivalence identifying the left adjoint of $\mathcal L^{\perp}\hookrightarrow \tcat$ with the Verdier quotient $\tcat\to \tcat/\mathcal L$.
\item[(e)] $\mathcal L=\Image(L)=\Ker(R)$ and $\mathcal R=\Ker(L)=\Image(R)$.
\qed
\end{itemize}
\end{prop}

\begin{defi} \label{defi:compl_pair}
Following \cite{meyernest-bc},
if $\mathcal L, \mathcal R\subset \tcat$ are thick subcategories satisfying condition ($*$) of Proposition \ref{prop:compl}, we say that $(\mathcal L,\mathcal R)$ is a pair of \emph{complementary subcategories of~$\tcat$}. The functorial distinguished triangle in ($*$) will be called the \emph{gluing triangle (at $A$)} for the complementary pair~$(\mathcal L,\mathcal R)$.
\end{defi}

We also recall the following immediate consequence of Proposition \ref{prop:compl}.

\begin{cor} \label{cor:compl_comp}
 If $(\mathcal L, \mathcal R)$ and $(\mathcal{\tilde L},  \mathcal{\tilde R})$ are two complementary pairs in $\tcat$ such that $\mathcal L\subset \mathcal{\tilde L}$ (equivalently: such that $\mathcal R\supset\mathcal{\tilde R}$) with gluing triangle $L\to \id\to R\to TL$, resp.\ $\tilde L\to \id \to \tilde R\to T\tilde L$, then 
$\tilde R\simeq \tilde RR$ and $L\tilde L\simeq L$.
\qed
\end{cor} 

One can use Brown representability to produce complementary pairs:

\begin{prop} \label{prop:brown_bousfield}
Let $\tcat$ be a triangulated category with small${}_{\alpha}$ coproducts. If $\mathcal S\subset \tcat_c$ is a small${}_{\alpha}$ subset of compact${}_{\alpha}$ objects, then $(\langle \mathcal S\rangle_{\loc} , \mathcal S^{\perp})$ is a complementary pair of localizing${}_{\alpha}$ subcategories of $\tcat$, depending only on the thick subcategory $\langle\mathcal S\rangle \subset \tcat_c$. \qed
\end{prop}

The proof of yet another well-known result, namely Neeman's localization theorem (\cite{neemanLoc}), also works verbatim in the $\alpha$-relative setting. 

\begin{thm}[Neeman localization theorem] \label{thm:neeman}
Let $\tcat$ be a compactly${}_{\alpha}$ generated triangulated category. Let $\mathcal L_0\subset \tcat_c$ be some (necessarily  essentially small${}_{\alpha}$) subset of compact${}_{\alpha}$ objects, and let $\mathcal L:=\langle\mathcal L_0\rangle_{\loc}$ be the localizing${}_{\alpha}$ subcategory of $\tcat$ generated by~$\mathcal L_0$.
Consider the resulting diagram of inclusions and quotient functors.
\begin{equation*}
\xymatrix{
\mathcal L \ar@{ >->}[r] & \tcat \ar@{->>}[r] & \tcat/\mathcal L \\
\mathcal L_c \ar@{ >->}[r] \ar@{ >->}[u] & \tcat_c \ar@{->>}[r] \ar@{ >->}[u] & \tcat_c/\mathcal L_c \ar@{..>}[u]_{F}
}
\end{equation*}
Then the following hold true:
\begin{itemize}
\item[\textrm{(a)}] The induced functor $F$ is fully faithful.
\item[\textrm{(b)}] The image of $F$ consists of compact${}_{\alpha}$ objects of $\tcat/\mathcal L$.
\item[\textrm{(c)}] $F(\tcat_c/\mathcal L_c)$ is a cofinal subcategory of $(\tcat/\mathcal L)_c$: for every $A\in (\tcat/\mathcal L)_c$ there are objects $A'\in(\tcat/\mathcal L)_c$ and $B\in \tcat_c/\mathcal L_c $ such that $A\oplus A'\simeq F(B)$.\qed
\end{itemize}
\end{thm}

Not everything generalizes, however. As the next example shows, arbitrary small${}_{\alpha}$ products are representable in a compactly${}_{\alpha}$ generated category only when $\alpha$ is inaccessible (which is, essentially, the case of a genuine compactly generated category). As a consequence, the representation theorem for \emph{covariant} functors (\cite[Thm.\ 2.1]{neemanDual}) is not available -- it cannot even be formulated in the usual way. See also Example \ref{ex:alpha} for a related problem.

\begin{example} \label{ex:products}
Let $\tcat$ be a compactly${}_{\alpha}$ generated triangulated category, and assume that the cardinal number $\alpha$ is \emph{not} inaccessible, i.e., that there exists a cardinal $\beta$ with $\beta<\alpha$ and $2^{\beta}\geq \alpha$ (e.g.\ $\alpha=\aleph_1$). If $0\not\simeq A\in \tcat_c$ is a nontrivial compact${}_{\alpha}$ object, then its $\beta$-fold product cannot exist in $\tcat$, because otherwise we would have $|\Hom(A,\prod_{\beta}A)| = |\prod_{\beta} \Hom (A,A)|\geq 2^{\beta} \geq \alpha$, in contradiction with the comp\-act${}_{\alpha}$\-ness of~$A$.
\end{example}

\subsection{The spectrum of a $\otimes$-triangulated category}
\label{subsec:spec}

We recall from \cite{balmer_prime} some basic definitions and results of Paul Balmer's geometric theory of tensor triangulated categories, or ``tensor triangular geometry''.

\begin{defi} \label{defi:tensor_tr}
By a \emph{tensor triangulated category} we always mean a triangulated category $\tcat$ (\cite{verdier} \cite{neemanTr}) equipped with a tensor product $\otimes: \tcat\times \tcat\to \tcat$ (i.e., a symmetric monoidal structure, see \cite{maclane}); we denote the unit object by~$\unit$. We assume that $\otimes$ is a triangulated functor in both variables, and we also assume that the natural switch $T(\unit)\otimes T(\unit)\stackrel{\sim}{\to}T(\unit)\otimes T(\unit)$ given by the tensor structure is equal to minus the identity.
Following \cite{balmer_spec2}, we call
\begin{eqnarray*}
\mathrm R_{\tcat}:=\End_{\tcat}(\unit) 
\quad \textrm{ and } \quad
\mathrm R_{\tcat}^*(\unit):=\End^*_{\tcat}(\unit):=\bigoplus_{n\in\Z}\Hom_{\tcat}(\unit,T^n(\unit))
\end{eqnarray*}
the \emph{central ring} and the \emph{graded central ring} of $\tcat=(\tcat,\otimes,\unit)$, respectively.
\end{defi}

\begin{remark} \label{rem:actions}
The central ring $\mathrm R_{\tcat}$ is commutative, and it acts on the whole category via $f\mapsto r\cdot f:= r\otimes f : A\simeq \unit \otimes A \to \unit \otimes B\simeq B$, for $r\in \mathrm R_{\tcat}$ and $f\in \Hom(A,B)$; we use here the structural identifications $\unit \otimes A\simeq A\simeq A\otimes \unit$. This makes $\tcat$ canonically into an $\mathrm R_{\tcat}$-linear category.
Our hypothesis on the switch $T(\unit)^{\otimes 2}\simeq T(\unit)^{\otimes 2}$ ensures that the graded central ring $\mathrm R_{\tcat}^*$ is graded commutative, by a classical argument. Also, it implies that the tensor product makes each graded Hom set $\Hom^*(A,B):= \bigoplus_n\Hom(A,T^nB)$ into a graded (left) module over $\mathrm R_{\tcat}^*$ such that composition is bilinear up to a sign rule (see \cite{balmer_spec2} or \cite[$\S$~2.1]{thesis} for details). In the following, we will localize these graded modules at homogeneous prime ideals $\mathfrak p$ of $\mathrm R^*_{\tcat}$, see~\ref{rem:gr_loc}.
\end{remark}

\begin{defi}[The spectrum] \label{defi:spectrum}
Let $\tcat$ be an essentially small  $\otimes$-triangulated category.
A \emph{prime tensor ideal} $\mathcal P$ in $\tcat$ is a proper (i.e.\ $\mathcal P \subsetneq \tcat$) thick subcategory of $\tcat$, which is a tensor ideal ($A\in \mathcal P,B\in \tcat$ $\Rightarrow$ $A\otimes B\in \mathcal P$) and is prime ($A\otimes B\in \mathcal P$ $\Rightarrow$ $A\in\mathcal P$ or $B\in \mathcal P$). The \emph{spectrum of $\tcat$}, denoted $\spc(\tcat)$, is the small set of its prime ideals.
The \emph{support} of an object $A\in \tcat$ is the subset
$$
\supp(A):=\{\mathcal P\mid A\not\in \mathcal P\}
 = \{\mathcal P\mid A\not\simeq 0 \textrm{ in } \tcat/\mathcal P \} \subset \spc(\tcat).
$$
We give the spectrum the \emph{Zariski topology}, which has $\{\spc(\tcat)\smallsetminus\supp(A)\}_{A\in \tcat}$ as a basis of open subsets. The space $\spc(\tcat)$ is naturally equipped with a sheaf of commutative rings $\Osheaf_{\mathcal T}$ whose stalks are the local rings $\Osheaf_{\mathcal T,\mathcal P}= \mathrm R_{\tcat/\mathcal P}$ (see \cite{balmer_spec2}). The resulting locally ringed space is denoted by $\spec(\tcat):=(\spc(\tcat),\Osheaf_{\tcat})$.
\end{defi}

\begin{remark} \label{rem:spectral}
The spectrum $\spc(\tcat)$ is a \emph{spectral space}, in the sense of Hochster \cite{hoch}: it is quasi-compact, its quasi-compact open subsets form an open basis, and every irreducible closed subset has a unique generic point. The support $A\mapsto \supp(A)$ is compatible with the tensor triangular structure, and is the finest such:
\end{remark}

\begin{prop}[Universal property {\cite{balmer_prime}}] \label{prop:UPsupp}
The support $A\mapsto \supp(A)$ has the following properties.
\begin{itemize}
\item[(SD1)] $\supp(0)=\emptyset$ and $\supp(\unit)=\spc(\tcat)$.
\item[(SD2)] $\supp(A\oplus B)= \supp(A)\cup \supp(B)$.
\item[(SD3)] $\supp(TA)=\supp(A)$.
\item[(SD4)] $\supp(B)\subset\supp(A)\cup \supp(C)$ if  $A\to B\to C\to TA$ is distinguished.
\item[(SD5)] $\supp(A\otimes B)=\supp(A)\cap\supp(B)$.
\end{itemize}
Moreover, if $(X,\sigma)$ is a pair consisting of a topological space $X$ together with an assignment $A\mapsto \sigma(A)$ from objects of $\tcat$ to closed subsets of $X$, satisfying the above five properties (in which case we say that $(X,\sigma)$ is a \emph{support datum} on $\tcat$), then there exists a unique \emph{morphism of support data} $f:(X,\sigma)\to (\spc(\tcat),\supp)$, i.e., a continuous map $f:X\to \spc(\tcat)$ such that $\sigma(A)=f^{-1}\supp(A)$ for all $A\in\tcat$. Concretely, $f$ is defined by $f(x):=\{A\in \tcat\mid x\not\in \sigma(A) \}$.
\qed
\end{prop}

\begin{terminology}
In the following, by ``a support'' $(X,\sigma)$ on some tensor triangulated category $\tcat$ we will simply mean a space $X$ together with some assignment $\sigma:\obj(\tcat)\to 2^X$ possibly lacking (some of) the good properties of a support datum.
\end{terminology}

Thus, the spectrum $(\spc(\tcat),\supp)$ is the universal support datum on~$\tcat$.  It has another important characterization. 

\begin{defi} \label{defi:thomason}
We say that a $\otimes$-ideal $\mathcal J\subset \tcat$ is \emph{radical} if $A^{\otimes n}\in \mathcal J$ for some $n\geq 1$ implies $A\in\mathcal J$. A subset $Y\subset \spc(\tcat)$ of the form $Y=\bigcup_i Z_i$, where each $Z_i$ is closed with quasi-compact open complement, is called a \emph{Thomason subset}.
\end{defi}

\begin{thm}[Classification theorem \cite{balmer_prime} \cite{bks}]
\label{thm:abstract_class}
The assignments 
\begin{eqnarray}
\mathcal J\mapsto \bigcup_{A\in\mathcal J} \supp(A)
&\textrm{ and } & 
Y \mapsto \{A\in \tcat\mid \supp(A)\subset Y\}
\end{eqnarray}
define mutually inverse bijections between the set of radical thick $\otimes$-ideals of $\tcat$ and the set of Thomason subsets of its spectrum $\spc(\tcat)$.

Conversely, if $(X,\sigma)$ is a support datum on $\tcat$ inducing the above bijection and  with $X$ spectral (in which case we say that $(X,\sigma)$ is a \emph{classifying} support datum), then the canonical morphism $f:(X,\sigma)\to (\spc(\tcat),\supp)$ is invertible; in particular, $f:X\to \spc(\tcat)$ is a homeomorphism. \qed
\end{thm}

So, up to canonical isomorphism, $(\spc(\tcat),\supp)$ is the unique classifying support datum on~$\tcat$. In examples so far, all explicit descriptions of the spectrum have been obtained from the Classification theorem, by proving that a specific concrete support datum is classifying.

\subsection{Rigid objects}

It often happens that the tensor product in a triangulated category is \emph{closed}, i.e., it has an internal Hom functor $\inthom:\tcat^{\op}\times \tcat\to \tcat$ providing a right adjoint $\inthom(A,?):\tcat\to\tcat$ of $?\otimes A:\tcat\to \tcat$ for each object $A\in \tcat$. 

 Being right adjoint to a triangulated functor, each $\inthom(A,?)$ is triangulated. Under some mild hypothesis, $\inthom$ preserves distinguished triangles also in the first variable: see \cite[App.\ C]{murfet} (I thank Amnon Neeman for the reference). In general, it is easily verified that the functor $\inthom(\,\textrm{?`},A)$ sends every distinguished triangle to a triangle that, while possibly not belonging to the triangulation, still yields long exact sequences upon application of the Hom functors $\Hom_{\tcat}(B,?)$.
The latter property suffices for many purposes, such as the proof of Prop.\ \ref{prop:thick_rigid} below.

\begin{example}
If $\tcat$ is a genuine compactly generated tensor triangulated category where the tensor commutes with coproducts, one obtains the internal Hom for free via Brown representability (simply represent the functors $\Hom_{\tcat}(\textrm{?`}\otimes A, B)$).
\end{example}

In the $\alpha$-relative setting, the internal Hom is only available when the source object is compact${}_{\alpha}$; fortunately, this suffices for our purposes. More precisely:

\begin{example} \label{ex:alpha}
Let $\tcat$ be a compactly${}_{\alpha}$ generated tensor triangulated category (Def.\ \ref{defi:cpt}) where $\otimes$ commutes with small${}_{\alpha}$ coproducts and where $\tcat_c\otimes \tcat_c\subset \tcat_c$. With these assumptions, if $A\in\tcat_c$ then Brown representability (Thm.\ \ref{thm:brown}) applies to the functor $\Hom(\textrm{?`}\otimes A,B):\tcat\to \Ab$, providing the right adjoint $\inthom(A,?):\tcat\to \tcat$ to tensoring with~$A$.
In general though there is a problem: if $\alpha$ is not inaccessible, i.e., if there exists a cardinal $\beta$ with $\beta<\alpha$ and $2^{\beta}\geq \alpha$ (e.g.\ $\alpha=\aleph_1$), then $\inthom$ \emph{cannot} be everywhere defined, as soon as $0\not\simeq\unit\in\tcat_{c}$. Indeed, if $X:=\inthom(\coprod_{\beta}\unit,\unit)\in\tcat$ were defined, we would have a natural isomorphism
$$
\Hom(A,X)
\simeq \Hom(A\otimes \coprod_{\beta}\unit,\unit)
\simeq \Hom(\coprod_{\beta}A, \unit)
\simeq \prod_{\beta}\Hom(A,\unit).
$$
Choosing $A=\unit \not\simeq 0$ we would obtain $|\Hom(\unit, X)|= |\prod_{\beta} \End (\unit)|\geq 2^{\beta}\geq \alpha$, contradicting the hypothesis that $\unit$ is compact${}_{\alpha}$. (Alternatively, we see that $X\simeq \prod_{\beta}\unit\in \tcat$, which is impossible by Example \ref{ex:products}).
\end{example}

\begin{defi} \label{defi:str_dual} Let $\tcat$ be a closed $\otimes$-triangulated category. We write $A^{\vee}:=\inthom(A,\unit)$ for the \emph{dual} of an object $A\in \tcat$.
An object $A\in \tcat$ is \emph{rigid} (or \emph{strongly dualizable}), if the morphism $A^{\vee}\otimes\, ? \to \inthom(A,?):\tcat\to \tcat$ -- obtained canonically by adjunction -- is an isomorphism. The $\otimes$-category $\tcat$ is \emph{rigid} if all its objects are rigid.
\end{defi}

\begin{prop}[{See \cite[App.\ A]{hps}}] \label{prop:thick_rigid}
Let $\tcat$ be a closed $\otimes$-triangulated category. The full subcategory of rigid objects is a thick $\otimes$-triangulated subcategory of $\tcat$ (in particular it contains the tensor unit). The contravariant functor $A\mapsto A^{\vee}$ restricts to a duality (i.e., $(?)^{\vee\vee}\simeq \id$) on this subcategory. \qed
\end{prop}

\begin{convention} \label{conv:cptly_gen}
We say that $\tcat=(\tcat,\otimes,\unit)$ is a \emph{compactly generated tensor triangulated category} if it is a tensor triangulated category (Def.\ \ref{defi:tensor_tr}) and if $\tcat$ is compactly${}_{\alpha}$ generated (Def.\ \ref{defi:cpt}) for some uncountable regular cardinal~$\alpha$, possibly with $\alpha=$ the cardinality of a proper class (what we dub the ``genuine'' case, that is, the usual sense of ``compactly generated''). Moreover, we assume that
\begin{itemize}
\item[(a)] for every $A\in \tcat$ the triangulated functors $A\otimes\,?$ and $?\otimes A$ preserve small${}_{\alpha}$ coproducts, and
\item[(b)] $\tcat_c\otimes \tcat_c\subset \tcat_c$ (cf.\ Ex.\ \ref{ex:alpha}) and the compact and rigid objects of $\tcat$ coincide.
\end{itemize}
In particular, $\tcat_c$ is a (rigid) tensor triangulated subcategory of $\tcat$.
\textbf{\emph{From now on, we will also drop the fixed cardinal $\alpha$ from our terminology.}}
\end{convention}

\begin{remark} \label{rem:genuine}
In the case of a genuine compactly generated category, as well as in the monogenic case (i.e., $\unit \in \tcat_c$ and $\tcat=\langle \unit \rangle_{\loc}$), the hypotesis $\tcat_c\otimes \tcat_c\subset \tcat_c$ is superfluous. Also, in general (and assuming (a)), to have equality of compact and rigid objects one needs only check that $\unit$ is compact and that $\tcat$ has a generating set consisting of compact and rigid objects.
\end{remark}

\begin{lemma}\label{lemma:ideal_ideal}
Let $\tcat$ be a compactly generated $\otimes$-triangulated category and $\mathcal J\subset \tcat_c$ a $\otimes$-ideal of its compact objects. Then $\langle \mathcal J\rangle_{\loc}$ is a localizing $\otimes$-ideal of~$\tcat$.
\end{lemma}

\proof
 For an object $A\in \tcat$, consider $\mathcal S_A:=\{X\in \tcat\mid X\otimes A\in \langle \mathcal J\rangle_{\loc}\}$. We must show that $\mathcal S_A=\tcat$ for all $A\in \langle \mathcal J\rangle_{\loc}$.
Note that $\mathcal S_A$ is always a localizing triangulated subcategory of~$\tcat$, because so is $\langle \mathcal J\rangle_{\loc}$ and because $\otimes$ preserves distinguished triangles and small coproducts. If $A\in \mathcal J$, then $\tcat_c\subset \mathcal S_A$ by hypothesis and therefore $\mathcal S_A= \tcat$. 
Now consider $\mathcal U:=\{A\in \tcat \mid \mathcal S_A=\tcat\}$. We have just seen that $\mathcal J\subset \mathcal U$, and one verifies immediately that $\mathcal U$ is a localizing subcategory of~$\tcat$. It follows that $\langle \mathcal J\rangle_{\loc}\subset \mathcal U$, as required.
\qed

The next result was first considered in stable homotopy by H.\ R.\ Miller \cite{miller}; cf.\ also \cite[Thm.\ 3.3.3]{hps} or \cite[Prop.\ 8.1]{bik}. In the topologist's jargon, it says that ``finite localizations are smashing''.
\begin{thm}[Miller] \label{thm:miller}
Let $\tcat$ be a compactly generated $\otimes$-triangulated category (as in Convention\ \ref{conv:cptly_gen}), and let $\mathcal J\subset \tcat_c$ be a tensor ideal of its compact objects. Then $\mathcal J^{\perp}=(\langle \mathcal J\rangle_{\loc})^{\perp}$ is a localizing tensor ideal, so that $(\langle \mathcal J\rangle_{\loc}, \mathcal J^{\perp})$ is a pair of complementary localizing tensor ideals of~$\tcat$.
\end{thm}

\proof
It follows from Prop.\ \ref{prop:brown_bousfield} that $(\langle \mathcal J\rangle_{\loc} , \mathcal J^{\perp})$ is a complementary pair of localizing subcategories, and from Lemma \ref{lemma:ideal_ideal} that $\langle \mathcal J\rangle_{\loc}$ is a $\otimes$-ideal of $\tcat$. It remains to see that $\mathcal J^{\perp}$ is a $\otimes$-ideal. Let $A\in \mathcal J^{\perp}$, and consider the full subcategory
$
\mathcal V_A:= \{X\in \tcat\mid X\otimes A \in \mathcal J^{\perp}  \}
$
of $\tcat$. It is triangulated and localizing because so is $\mathcal J^{\perp}$. It contains every compact object: if $C\in \tcat_c$ and $J\in \mathcal J$, then $\Hom(J,C\otimes A)\simeq \Hom(J\otimes C^{\vee},A)\simeq 0$  because $C$ is rigid and $\mathcal J$ is an ideal. Therefore $\mathcal V_A= \langle \tcat_c \rangle_{\loc}= \tcat$, that is to say $\tcat\otimes A\subset \mathcal J^{\perp}$, for all $A\in \mathcal J^{\perp}$.
\qed

\begin{remark} \label{remark:tensor_gluing}
If both subcategories $\mathcal L,\mathcal R\subset \tcat$ in a complementary pair $( \mathcal L, \mathcal R)$ are $\otimes$-ideals, then the gluing triangle for an arbitrary object $A\in \tcat$ is obtained by tensoring $A$ with the gluing triangle for the $\otimes$-unit~$\unit$. (This is an exercise application of the uniqueness of the gluing triangle, see Prop.\ \ref{prop:compl}.) 
\end{remark}

\subsection{Central localization}
\label{subsec:central}

In a tensor triangulated category $\tcat$, as we already mentioned, the tensor product naturally endows the Hom sets with an action of the central ring $\mathrm R_{\tcat}=\End_{\tcat}(\unit)$, making $\tcat$ an $\mathrm R_{\tcat}$-linear category. If $S\subset \mathrm R_{\tcat}$ is a multiplicative system, one may localize each Hom set at $S$. As the next theorem shows, the resulting category still carries a tensor triangulated structure. Let us be more precise.

\begin{construction} \label{defi:gen_loc}
Let $\mathcal C$ be an $R$-linear category, for some commutative ring~$R$.
Let $S\subset R$ be a multiplicative system  (i.e., $1\in S$ and $S\cdot S\subset S$). Define $S^{-1}\mathcal C$ to be the category with the same objects as $\mathcal C$, with Hom sets the localized modules $S^{-1}\Hom_{\mathcal C}(A,B)$ and with composition defined by $ (\frac{g}{t} , \frac{f}{s} )\mapsto \frac{g\circ f}{ts}$.
One verifies easily that $S^{-1}\mathcal C$ is an $S^{-1}R$-linear category and that there is an $R$-linear canonical functor $\loc : \mathcal C\to S^{-1}\mathcal C$. It is the universal functor from $\mathcal C$ to an $S^{-1}R$-linear category.
\end{construction}

\begin{defi} \label{defi:central_loc}
Let $\tcat$ be a tensor triangulated category, and let $S\subset \mathrm R_{\tcat}$ be a multiplicative system of its central ring. We call $S^{-1}\tcat$ (as in \ref{defi:gen_loc}) the \emph{central localization of $\tcat$ at~$S$}. The next result shows that it is again a tensor triangulated category.
\end{defi}

\begin{thm}[Central localization {\cite[Thm.\ 3.6]{balmer_spec2}}] \label{thm:central_loc}
Consider the thick $\otimes$-ideal $\mathcal J=\langle \cone(s)\mid s\in S \rangle_{\otimes}\subset \tcat$ generated by the cones of maps in~$S$. Then there is a canonical isomorphism $S^{-1}\tcat\simeq \tcat/\mathcal J$ which identifies $\loc: \tcat\to S^{-1}\tcat$ with the Verdier quotient $q:\tcat\to \tcat/\mathcal J$. In particular, the central localization $S^{-1}\tcat$ inherits a canonical $\otimes$-triangulated structure such that $\loc$ is $\otimes$-triangulated; conversely, $q$ is the universal $R$-linear triangulated functor to an $S^{-1}R$-linear $\otimes$-triangulated category.
\qed
\end{thm}

The procedure of central localization can be adapted to compactly generated categories in a most satisfying way, as we expound in the next theorem.

\begin{thm} \label{thm:lochom}
Let $\tcat$ be a compactly generated $\otimes$-triangulated category (as in \ref{conv:cptly_gen}), and let $S$ be a multiplicative subset of the central ring~$\mathrm R_{\tcat}$. Write 
\begin{equation*}
\mathcal J_S:=\langle\cone(s)\mid s\in S \rangle_{\otimes} \subset \tcat_c\quad,\quad
\mathcal L_S:=\langle \mathcal J_S\rangle_{\loc} \subset \tcat  .
\end{equation*}
The objects of $\mathcal T_S:=(\mathcal L_S)^{\perp}$ will be called \emph{$S$-local objects}. Then the pair 
$ 
( \mathcal L_S \;, \; \tcat_S )
$ 
is a complementary pair (Def.~\ref{defi:compl_pair}) of localizing $\otimes$-ideals of $\tcat$. In particular, the gluing triangle for an object $A\in\tcat$ is obtained by tensoring $A$ with the gluing triangle for the $\otimes$-unit
\begin{equation*} 
\xymatrix@1{
L_S(\unit) \ar[r]^-{\varepsilon} & \unit \ar[r]^-{\eta} & R_S(\unit) \ar[r] & TL_S(\unit) .
}
\end{equation*}
This situation has the following properties:
\begin{itemize}
\item[(a)] $\mathcal L_S=L_S(\unit)\otimes \tcat$ and  $\mathcal T_S= R_S(\unit) \otimes \tcat$.
\item[(b)] $\varepsilon:L_S(\unit)\otimes L_S(\unit) \simeq L_S(\unit) $ and $\eta:R_S(\unit)\simeq R_S(\unit) \otimes R_S(\unit)$.
\item[(c)] $\tcat_S$ is again a compactly generated $\otimes$-triangulated category, as in Conv.\ \ref{conv:cptly_gen}, with tensor unit $R_S(\unit)$. (Note that $R_S(\unit)$ is compact in $\tcat_S$, but need not be in $\tcat$.)
\item[(d)] Its compact objects are $(\tcat_S)_c=\langle R_S(\tcat_c)\rangle\subset \tcat_S$. (Again, they are possibly non compact in $\tcat$.)
\item[(e)] The functor 
 $ R_S= R_S(\unit)\otimes\,?:\tcat\to \mathcal \tcat_S $ 
is an $\mathrm R_{\tcat}$-linear $\otimes$-triangulated functor commuting with small coproducts. It takes generating sets to generating sets.
\item[(f)] To apply $\Hom(\unit,?)$ on $\unit\stackrel{\eta}{\to} R_S(\unit)$ induces the localization $\mathrm R_{\tcat} \to S^{-1}\mathrm R_{\tcat}$. It follows in particular that $\mathrm R_{\tcat_S}= S^{-1}\mathrm R_{\tcat}$.
\item[(g)] An object $A\in \tcat$ is $S$-local if and only if $s\cdot \id_A$ is invertible for every $s\in S$.
\item[(h)] If $A\in \tcat_c$, then $\eta:B\to R_S(\unit)\otimes B$ induces an isomorphism
\begin{equation*}
S^{-1}\Hom_{\tcat}(A,B)\simeq \Hom_{\tcat}(A,R_S(\unit)\otimes B)
\end{equation*}
for every $B\in \tcat$.
\end{itemize}
\end{thm}

\begin{remarks}
(a)
The category $\mathcal L_S$ is both compactly generated and a tensor triangulated category but, since in general its $\otimes$-unit $L_S(\unit)$ is not compact, it may fail to be a compactly generated tensor triangulated category as defined in Convention \ref{conv:cptly_gen}.

(b)
There are graded versions of the above results, where one considers multiplicative systems of the graded central ring $\mathrm R^*_{\tcat}=\End^*(\unit)$. We don't use them here, so we have omitted their (slightly more complicated) formulation.

(c) We don't really need that all compact objects be rigid (as was assumed in Convention \ref{conv:cptly_gen}) in order to prove Theorem \ref{thm:lochom}. More precisely, one can show that $\tcat_S$ is a $\otimes$-ideal in $\tcat$ without appealing to Miller's Theorem. It suffices to use the $\mathrm R_{\tcat}$-linearity of the tensor product and the characterization of $S$-local objects (part (g) of the theorem): if $A\in \tcat_S$ and $B\in \tcat$, then $s\cdot \id_{A\otimes B}= (s\cdot \id_A) \otimes B$ is invertible for all $s\in S$ and therefore $A\otimes B\in \tcat_S$.
\end{remarks}

\proof[Proof of Theorem \ref{thm:lochom}]

The first claim is Miller's Theorem \ref{thm:miller} and Remark \ref{remark:tensor_gluing}, applied to the $\otimes$-ideal $\mathcal J_S\subset \tcat_c$. Thus $(\mathcal L_S, \tcat_S)$ is a complementary pair of localizing $\otimes$-ideals. Part (a) and (b) are then formal consequences.  The statements in (c)-(e) are either clear, or follow from Neeman's Localization Theorem \ref{thm:neeman} (the $\mathrm R_{\tcat}$-linearity in (e) is Lemma \ref{lemma:action-LR} below). Let's now turn to the more specific claims (f)-(h).

\begin{lemma} \label{lemma:q-inverts-S}
The quotient functor $q:\tcat\to \tcat/\mathcal L_S$ is $\mathrm R_{\tcat}$-linear and it inverts all endomorphisms of the form $s\cdot\id_A$ with $s\in S$ and $A\in\tcat$.
\end{lemma}

\proof
Let $s\in S$ and $A\in\tcat$. Then $\cone(s\cdot\id_A)= \cone(s)\otimes A$ belongs to $\mathcal L_S $, because $\cone(s)\in\mathcal J_S\subset \mathcal L_S$ by definition and $\mathcal L_S$ is a $\otimes$-ideal. 
\qed

In particular, by the universal property of central localization \eqref{defi:gen_loc}, the quotient functor $q:\tcat\to\tcat/ \mathcal L_S$ factors as
\begin{equation*}
\xymatrix{\tcat \ar[r]^-{q}  \ar[d]_{\loc}  & \tcat/\mathcal L_S \\
S^{-1}\tcat .\ar@{..>}[ur]_{\overline{q}} &
}
\end{equation*}
We clearly have a commutative square
\begin{equation} \label{centrloc-neemanloc}
\xymatrix{
S^{-1}\tcat \ar[r]^-{\overline{q}} & \tcat/{\mathcal L_S} \\
S^{-1}\tcat_c \ar@{ >->}[u] \ar[r]^-{\overline{q}_c}_{\simeq} & \tcat_c/\mathcal J_S \ar@{ >->}[u]
}
\end{equation}
where every functor is the identity or an inclusion on objects, and where $\overline q_c$ is the canonical identification of Theorem \ref{thm:central_loc}; the right vertical functor is fully faithful by Theorem \ref{thm:neeman}~(a).

\begin{prop} \label{prop:mythm47}
The canonical functor $\overline{q}$ restricts to an isomorphism
\begin{equation*}
\overline{q} \;:\; S^{-1}\Hom_{\tcat}(C,B)\stackrel{\sim}{\longrightarrow} \Hom_{\tcat/\mathcal L_S}(C,B)
\end{equation*}
of $S^{-1}\mathrm R_{\tcat}$-modules for all compact $C\in \tcat_c$ and arbitrary $B\in\tcat$. 
\end{prop}

\proof
Fix a $C\in\tcat_c$. We may view
\begin{equation} \label{ellC}
\overline{q}:S^{-1}\Hom_{\tcat}(C,?)\longrightarrow \Hom_{\tcat/\mathcal L_S}(C,?)
\end{equation}
as a morphism of homological functors to $S^{-1}\mathrm R_{\tcat}$-modules, both of which commute with small coproducts. Moreover, $\overline{q}$ is an isomorphism on compact objects, as we see from \eqref{centrloc-neemanloc}.
 It follows that \eqref{ellC} is an isomorphism on the localizing subcategory generated by $\tcat_c$, which is equal to the whole category~$\tcat$. 
\qed

Part (h) of the theorem is now an easy consequence, provided we correctly identify the isomorphism in question.

\begin{cor} \label{cor:varloc}
Let $C,B\in\tcat$ with $C$ compact. Then $\eta_B:B\to R_S(B)$ induces an isomorphism
$ 
\beta\; : \; S^{-1}\Hom_{\tcat}(C,B) \stackrel{\sim}{\longrightarrow} \Hom_{\tcat}(C,R_S(B))
$ 
of $\mathrm R_{\tcat}$-modules. 
\end{cor}

\proof 
Recall from \ref{prop:compl} (c)-(d) that $q$ has a fully faithful right adjoint $q_r$ such that $R_S= q_rq$. 
Since $\eta$ is natural, the following square commutes for all $f:C\to B$,
\begin{equation*}
\xymatrix{
C\ar[d]_f \ar[r]^-{\eta_C} & q_rq(C) \ar[d]^{q_rq(f)} \\
B \ar[r]^-{\eta_B} & q_rq(B)
}
\end{equation*}
showing that the next (solid) square is commutative. 

\begin{equation*}
\xymatrix{
& \Hom_{\tcat}(C,B) \ar[dl]_{\loc} \ar[dd]^{(\eta_B)_*} \ar[r]^-{q} & \Hom_{\tcat/\mathcal L_S}(qC,qB) \ar[dd]^{q_r}_{\simeq} \\
S^{-1}\Hom_{\tcat}(C,B)\ar@{..>}[urr]_<<<<<<{\overline{q}} \ar@{..>}[dr]_{\beta} && \\
& \Hom_{\tcat}(C,q_rqB) & \Hom_{\tcat}(q_rqC,q_rqB) \ar[l]_-{(\eta_C)^*}^-{\simeq}
}
\end{equation*}
Notice that $(\eta_C)^*$ is an isomorphism by \ref{prop:compl}~(a). By the compactness of $C$ and by Proposition~\ref{prop:mythm47}, $q$ induces the isomorphism~$\overline{q}$. Composing this isomorphism with the other two, we see that~$\beta$, the factorization of $(\eta_B)_*$ through $\loc$, is an isomorphism as claimed.
\qed

\begin{lemma} \label{lemma:action-LR}
The endofunctors $L_S$ and $R_S$ are $\mathrm R_{\tcat}$-linear.
\end{lemma}

\proof
This can be seen in various ways. For instance, by applying the functorial gluing triangle $L_S\to \id \to R_S \to TL_S$ to $r\cdot f :A\to B$, resp.\ by applying it to $f:A\to B$ and then multiplying by~$r$, we obtain two commutative squares
\begin{equation*}
\xymatrix@1{
A \ar[d]_{r\cdot f} \ar[r]^-{\eta_A} & R_SA \ar[d]^{R_S(r\cdot f)} \\
B \ar[r]^-{\eta_B} & R_SB
}
\quad \quad
\xymatrix@1{
A \ar[d]_{r\cdot f} \ar[r]^-{\eta_A} & R_SA \ar[d]^{r\cdot R_S( f)} \\
B \ar[r]^-{\eta_B} & R_SB.
}
\end{equation*} 
In particular, we see that the difference $d:=R_S( r\cdot  f) - r\cdot R_S( f)$ composed with $\eta_A$ is zero, so it must factor through $TL_SA\in \mathcal L_S$. But the only map $TL_SA\to R_SB$ is zero, hence $d=0$, that is $R_S( r\cdot  f) = r\cdot R_S( f)$. A similar argument applies to show that $L_S$ is $\mathrm R_{\tcat}$-linear.
\qed

Together with Lemma \ref{lemma:q-inverts-S}, the next lemma provides part~(g).

\begin{lemma} \label{lemma:converseS}
If $A\in \tcat$ is such that $s\cdot\id_A$ is invertible for all $s\in S$, then $\eta_A:A\to R_S(A)$ is an isomorphism. In particular, $A\in \Image(R_S)= \tcat_S$.
\end{lemma}

\proof
The map $\eta_A:A\to R_SA$ induces the following commutative diagram of natural transformations between cohomological functors $\tcat^{\op}\to \mathrm R_{\tcat}\Modules$:
\begin{equation*}
\xymatrix{
\Hom_{\tcat}(\textrm{?`},A) \ar[dr]_{\loc} \ar[rr]^-{(\eta_A)_*} && \Hom_{\tcat}(\textrm{?`},R_SA) \\
&  S^{-1}\Hom_{\tcat}(\textrm{?`},A) \ar[ur]_{\beta}  &
}
\end{equation*}
The hypothesis on $A$ implies that $\loc$
is an isomorphism. By Corollary~\ref{cor:varloc}, the map
$\beta$ is an isomorphism on compact objects. Hence their composition $(\eta_A)_*$ is a morphism of cohomological functors both of which send coproducts to products -- indeed they are representable -- and such that it is an isomorphism at each $C\in\tcat_c$. 
It follows that  $(\eta_A)_*$ is an isomorphism at every object. By Yoneda, $\eta_A$ is an isomorphism in~$\tcat$, showing that $A\in \Image(R_S)$.
\qed

Finally, part (f) is (h) for $A=B=\unit$; note for the second assertion that $\Hom(\unit, R_S(\unit)) $ $\simeq$ $ \Hom(R_S(\unit),R_S(\unit))=\mathrm R_{S^{-1}\tcat}$.
This ends the proof of Theorem \ref{thm:lochom}.
\qed

\begin{remark} 
The authors of \cite{bik} prove very similar results (and much more) for genuine compactly generated categories, without need for a tensor structure. Instead of the central ring $\mathrm R_{\tcat}$, they posit a noetherian graded commutative ring acting on $\tcat$ via endomorphisms of $\id_{\tcat}$, compatibly with the translation. If $\tcat$ is moreover a \emph{tensor} triangulated category (with our same hypotheses \ref{conv:cptly_gen}), they also prove the results in Theorem~\ref{thm:lochom} for the graded central ring $\mathrm R^*_{\tcat}$, but only when the latter is noetherian; see \cite[$\S$8]{bik}). Wishing to apply their results, we met the apparently insurmountable problem that in the $\alpha$-relative setting Brown representability for the dual, which is crucially used in \emph{loc.\ cit.}, is not available (cf.\ Ex.\ \ref{ex:products}). 
\end{remark}

\section{Classification in compactly generated categories} 
\label{sec:criterium}

\subsection{An abstract criterion}
Let $\mathcal K$ be an essentially small $\otimes$-triangulated category. In most examples so far where the Balmer spectrum $\spc(\mathcal K)$ has been described explicitly, $\mathcal K$ is  the subcategory $\tcat_c$ of compact and rigid objects in some compactly generated $\otimes$-triangulated category~$\tcat$. Indeed, the ambient category $\tcat$ provides each time essential tools for the computation of $\spc(\mathcal K)$. The next theorem, abstracted from the example of modular representation theory (see Example \ref{ex:class_bcr}), yields a general method for precisely this situation. 

\begin{thm} \label{thm:abstract}
Let $\tcat$ be a compactly generated $\otimes$-tri\-angulat\-ed category, as in Convention \ref{conv:cptly_gen}. Let $X$ be a spectral topological space, and let $\sigma: \obj(\tcat)\to 2^X$ be a function assigning to every object of $\tcat$ a subset of~$X$. Assume that the pair $(X,\sigma)$ satisfies the following ten axioms:
\begin{itemize}
\item[(S0)] $\sigma(0)=\emptyset$.
\item[(S1)] $\sigma(\unit)=X$.
\item[(S2)] $\sigma(A\oplus B)= \sigma(A)\cup \sigma(B)$.
\item[(S3)] $\sigma(TA)=\sigma(A)$.
\item[(S4)] $\sigma(B)\subset\sigma(A)\cup \sigma(C)$ for every distinguished triangle  $A\to B\to C\to TA$.
\item[(S5)] $\sigma(A\otimes B)=\sigma(A)\cap\sigma(B)$ for every compact $A\in \tcat_c$ and arbitrary $B\in \tcat$. 
\item[(S6)] $\sigma(\coprod_i A_i)= \bigcup_i \sigma(A_i)$ for every small family $\{A_i\}_i\subset \tcat$.
\item[(S7)] $\sigma(A)$ is closed in $X$ with quasi-compact complement $X\smallsetminus \sigma(A)$ for all  $A\in \tcat_c$.
\item[(S8)] For every closed subset $Z\subset X$ with quasi-compact complement, there exists an $A\in \tcat_c$ with $\sigma(A)=Z$.
\item[(S9)] $\sigma(A)=\emptyset \Rightarrow A\simeq 0$.
\end{itemize}
Then the restriction of $(X,\sigma)$ to $\tcat_c$ is a classifying support datum, so that, by Theorem \ref{thm:abstract_class}, the induced canonical map $X\to \spc(\tcat_c)$ is a homeomorphism.
\end{thm}

\begin{example} \label{ex:class_bcr}
Let $G$ be a finite group and $k$ a field. Let $\tcat$ be the stable module category $\stmod(kG):=\mathsf{mod}(kG)/\mathsf{proj}(kG)$ of finitely generated $kG$-modules, equipped with the tensor product $\otimes:=\otimes_k$ (with diagonal $G$-action) and the unit object $\unit:=k$ (with trivial $G$-action); see \cite{carlson}. Then there is a homeomorphism $\spc(\stmod(kG))\simeq \proj(H^*(G;k))$. 

Indeed, we may embed $\stmod(kG)$ as the full subcategory of compact and rigid objects inside $\StMod(kG)$, the stable category of possibly infinite dimensional $kG$-modules. The latter is a (genuine) compactly generated category as in \ref{conv:cptly_gen}; cf.\ e.g.\ \cite{rickard} \cite[$\S$10]{bik}.
Let $R:=H^*(G;k)=\End^{\geq0}_{\stmod(kG)}(k,k)$ be the cohomology ring of~$G$. 
Let $X:=\proj(H^*(G;k))
= \specgr(H^*(G;k))\smallsetminus \{\mathfrak m \}$, where $\mathfrak m=H^{>0}(G;k)$.
Consider on $\StMod(kG)$ the support $\sigma: \obj(\tcat)\to 2^X$ given by the \emph{support variety} of a module $M\in \StMod(kG)$, as introduced in \cite{bcrII}. It follows from the results of \emph{loc.\ cit.} that $(X,\sigma)$ satisfies all of our axioms (S0)-(S9). Most non-trivially, (S5) holds by the Tensor Product theorem \cite[Thm.\ 10.8]{bcrII} and (S9) by, essentially, Chouinard's theorem. Therefore by Theorem \ref{thm:abstract} there is a unique isomorphism $(X,\sigma) \simeq (\spc(\stmod(kG)), \supp)$ of support data on $\stmod(kG)$.


\end{example}

Before we give the proof of the theorem, we note that a common way of obtaining supports $(X,\sigma)$ on $\tcat$ is by constructing a suitable family of homological functors $F_x: \mathcal T\to \mathcal A_x$, $x\in X$. We make this intuition precise in the following -- somewhat pedant -- lemma, whose proof is a series of trivial verifications left to the reader.

\begin{lemma} \label{lemma:hml_supp}
Consider a family $\mathcal F=\{F_x:\tcat\to \mathcal A_x \}_{x\in X}$ of functors parametrized by a topological space~$X$. Assume that each $\mathcal A_x$ has a zero object $0$ (i.e., $0$ is initial and final in $\mathcal A_x$). For each $A\in \tcat$ we define
\begin{equation*}
\sigma_{\mathcal F}(A) := \{ x\in X \mid F_x(A)\not \simeq 0 \textrm{ in } \mathcal A_x\}\quad \subset \quad X.
\end{equation*}
Then, if the functors $\mathcal F=\{F_x\}_x$ satisfy condition \emph{(F$n$)} of the following list, the induced support  $(X,\sigma_{\mathcal F})$ satisfies the corresponding hypothesis \emph{(S$n$)} of Theorem~\ref{thm:abstract}.
\begin{itemize}
\item[(F0)] $F_x(0)\simeq 0\in \mathcal A_x$.
\item[(F1)] $F_x(\unit)\not\simeq 0 \in \mathcal A_x$.
\item[(F2)] $\mathcal A_x$ is additive and $F_x$ is an additive functor (thus (F2) $\Rightarrow$ (F0)).
\item[(F3)] $\mathcal A_x$ is equipped with an endoequivalence $T$ and $F_xT\simeq TF_x$. 
\item[(F4)] $\mathcal A_x$ is abelian and $F_xA\to F_xB\to F_xC$ is exact for every distinguished triangle $A\to B\to C\to TA$.
\item[(F5)] $\mathcal A_x =(\mathcal A_x,\,\hat\otimes\,)$ is a tensor category such that
\begin{eqnarray*}
M\,\hat\otimes\, N\simeq 0 &\Leftrightarrow& M\simeq 0 \textrm{ or } N\simeq 0 ,
\end{eqnarray*}
and there exist isomorphisms 
\begin{eqnarray*}
F_x(A \otimes B) &\simeq & F_x(A)\,\hat\otimes\, F_x(B) 
\end{eqnarray*}
at least for $A\in \tcat_c$ compact and $B\in \tcat$ arbitrary.
\item[(F6)] $F_x$ preserves small coproducts.
\item[(F9)] The family $\mathcal F=\{F_x\}_{x\in X}$ detects objects, i.e.: $F_x(A)\simeq 0\;\forall x \Rightarrow A\simeq 0$. \qed
\end{itemize}
\end{lemma}

A functor $F$ with properties (F2), (F3) and (F4) is usually called a \emph{stable homological functor} (also recalled in Def.\ \ref{defi:stable} below). Note also that the only \emph{collective} property of the family $\mathcal F$ is (F9).

In this generality, the translations of conditions (S7) and (S8) remain virtually identical, so we omitted them from our list (but see Prop.\  \ref{prop:S7S8} below for the discussion of a significant special case).

Let us now prove Theorem \ref{thm:abstract}. For any subset $Y\subset X$, let us use the notation 
\begin{eqnarray*}
\mathcal C_Y &:= & \{A\in \tcat_c \mid \sigma(A)\subset Y \} \;\; \subset \;\; \tcat_c \\
\tcat_Y&:=& \langle \mathcal C_Y \rangle_{\loc} \quad \subset \quad \tcat.
\end{eqnarray*}
We begin with some easy observations:

\begin{lemma} \label{lemma:suppZ}
\begin{itemize}
\item[(a)]
The subcategory $\mathcal C_Y\subset \tcat_c$ is a radical thick $\otimes$-ideal. In particular, it is a thick triangulated subcategory and thus $\mathcal C_Y=(\tcat_Y)_c$.
\item[(b)] 
If $A\in \mathcal T_Y$, then $\sigma(A)\subset Y$.
\end{itemize}
\end{lemma} 
\proof
(a) It follows immediatly from axioms (S0) and  (S2)-(S5) that $\mathcal C_Y$ is a thick triangulated tensor ideal of $\tcat_c$. Now let $A\in \tcat_c$ with $A^{\otimes n}\in \mathcal C_Y$ for some $n\geq 1$. This means $\sigma(A^{\otimes n})\subset Y$ and therefore $\sigma(A)\subset Y$ by (S5). Thus $\mathcal C_Y$ is radical.

(b) By the axioms (S0), (S2)-(S4) and (S6), the full subcategory $\{A\in \tcat\mid \sigma(A)\subset Y \}$ of all objects supported on $Y$ is a localizing triangulated subcategory of~$\tcat$. Since it obviously contains $\mathcal C_Y$, it must contain $\tcat_Y=\langle \mathcal C_Y\rangle_{\loc}$.
\qed

\begin{lemma}[{cf.\,\cite[Prop.\ 3.3]{bcr}}] \label{lemma:keycat} 
Let $\mathcal E\subset \tcat_c$ be any self-dual collection of compact objects, meaning that $\mathcal E=\mathcal E^{\vee}:=\{E^{\vee}\mid E\in \mathcal E\}$, and let $\sigma(\mathcal E):=\bigcup_{E\in \mathcal E}\sigma(E)\subset X$ denote their collective support. Then
$$
\langle \mathcal E\rangle_{\otimes}=\mathcal C_{\sigma(\mathcal E)}
$$
in $\tcat_c$, that is, the thick $\otimes $-ideal of $\tcat_c$ generated by $\mathcal E$ consists precisely of the compact objects which are supported on $\sigma(\mathcal E)$. 
\end{lemma}
 
\proof
Let us write $Y:= \sigma(\mathcal E)$.
Each of the thick subcategories $\langle \mathcal E\rangle_{\otimes}$ and $\mathcal C_Y $ of $ \tcat_c$ determines a complementary pair  in $\tcat$ by Proposition \ref{prop:brown_bousfield}, namely $(\langle \mathcal E \rangle_{\otimes,\loc}, \langle \mathcal E \rangle_{\otimes,\loc}^{\perp})$ and $(\tcat_Y,\tcat_Y^{\perp})$, with gluing triangles
\begin{eqnarray*}
\xymatrix@1{
L_{\langle \mathcal E\rangle_{\otimes}} \ar[r]
 & \id_{\tcat} \ar[r]
  & R_{\langle \mathcal E\rangle_{\otimes}} \ar[r]
   & TL_{\langle \mathcal E\rangle_{\otimes}}
}
&& \textrm{and}\\
\xymatrix@1{
L_{\mathcal C_Y} \ar[r]
 & \id_{\tcat} \ar[r]
  & R_{\mathcal C_Y} \ar[r]
   & TL_{\mathcal C_Y}
}
&&,
\end{eqnarray*}
respectively. Moreover, the two thick subcategories can be recovered as
\begin{equation*} \label{Ims}
\langle \mathcal E\rangle_{\otimes}=( \Image(L_{\langle \mathcal E\rangle_{\otimes}}) )_c \quad \textrm{and}\quad \mathcal C_Y =( \Image(L_{\mathcal C_Y}))_c. 
\end{equation*}
Thus, in order to prove the lemma, it suffices to find an isomorphism $L_{\langle\mathcal E\rangle_{\otimes}}\simeq L_{\mathcal C_Y}$.
Since $\mathcal C_Y$ is a thick $\otimes$-ideal (by Lemma \ref{lemma:suppZ}~(a)) and it contains $\mathcal E$, we must have the inclusion $\langle \mathcal E\rangle_{\otimes}\subset \mathcal C_Y$ and thus $\langle\mathcal E\rangle_{\otimes,\loc}\subset \tcat_{Y}$. It follows from Corollary \ref{cor:compl_comp} that $L_{\langle \mathcal E\rangle_{\otimes}} L_{\mathcal C_Y}\simeq L_{\langle \mathcal E\rangle_{\otimes}}$.
Hence, for any $A\in \tcat$, the first of the above gluing triangles applied to the object $L_{\mathcal C_Y}(A)$ takes the form
 \begin{equation} \label{gtr} 
\xymatrix@1{
L_{\langle \mathcal E \rangle_{\otimes}}(A) \ar[r]^{} & L_{\mathcal C_Y}(A) \ar[r]^{} & R_{\langle \mathcal E \rangle_{\otimes}}L_{\mathcal C_Y}(A) \ar[r] & TL_{\langle \mathcal E \rangle_{\otimes}}(A).
}
\end{equation}
Since $A\in \tcat$ is arbitrary, we have reduced the problem to proving that the third object $B:=R_{\langle \mathcal E \rangle_{\otimes}}L_{\mathcal C_Y}(A)$ in the distinguished triangle \eqref{gtr} is zero. By axiom (S9), it suffices to prove the following
\begin{itemize}
\item[\emph{Claim:}]  $\sigma(B)=\emptyset$.
\end{itemize}
Indeed, since the first two objects in \eqref{gtr} belong to the triangulated category~$\mathcal T_Y$, so does~$B$. Therefore $\sigma(B)\subset Y$ by Lemma \ref{lemma:suppZ}~(b). Let $E\in \mathcal E$, and let $C$ be any compact object of~$\mathcal T$.
Then
\begin{eqnarray*}
\Hom_{\mathcal T} (C,E^{\vee} \otimes B) \simeq \Hom_{\mathcal T}(C\otimes E, B)\simeq 0
\end{eqnarray*}
because $E\in \tcat_c$ is rigid (for the first isomorphism), and because $C\otimes E\in \langle \mathcal E\rangle_{\otimes}$ and $B\in \Image(R_{\langle \mathcal E\rangle_{\otimes}})= \langle \mathcal E\rangle_{\otimes}^{\perp}$ (for the second one). But this implies $E^{\vee}\otimes B\simeq 0$, because compact objects generate~$\mathcal T$. Hence $\sigma(E^\vee \otimes B)=\emptyset$ by~(S0). Using this fact, together with $\sigma(B)\subset Y=\sigma(\mathcal E)= \sigma(\mathcal E^{\vee})$, we conclude that
\begin{equation*}
\sigma(B)
= 
\left(\bigcup_{E\in \mathcal E}\sigma(E^{\vee})\right)\cap\sigma(B)
\stackrel{}{=}
\bigcup_{E\in \mathcal E} \sigma(E^{\vee})\cap \sigma(B)
\stackrel{\textrm{(S5})}{=}
\bigcup_{E\in \mathcal E} \sigma(E^{\vee}\otimes B)
=
\emptyset
\end{equation*}
as we had claimed.
\qed

\begin{lemma} \label{lemma:self_dual}
Every thick $\otimes$-ideal of $\tcat_c$ is self-dual.
\end{lemma}

\proof
This is \cite[Prop.\ 2.6]{balmer_filtr}; note that the hypothesis in \emph{loc.\ cit.\ }that the duality functor $(\cdot)^{\vee}$ be triangulated is not used in the proof. Indeed, let
 $\mathcal C\subset \tcat_c$ be a thick $\otimes$-ideal. Every rigid object $A$ is a retract of $A\otimes A^{\vee}\otimes A$ (this holds in any closed tensor category, by one of the triangular identities of the adjunction between $?\otimes A$ and $A^{\vee}\otimes?$
).
Then also $A^{\vee}$ is a direct summand of $A^{\vee}\otimes A^{\vee\vee}\otimes A^{\vee} \simeq A^{\vee}\otimes A\otimes A^{\vee}$. Since $\mathcal C$ is thick and $(\cdot)^{\vee}:\tcat_c\to \tcat^{\op}_c$ is an additive tensor equivalence, both $\mathcal C$ and $\mathcal C^{\vee}$ are closed under taking summands and tensoring with arbitrary objects of~$\tcat_c$. It follows from the previous remarks that $\mathcal C\subset \mathcal C^{\vee}$ and $\mathcal C^{\vee} \subset \mathcal C$.
\qed

\proof[Proof of Theorem \ref{thm:abstract}] By properties (S0)-(S5) and (S7), the restriction of $(X,\sigma)$ to $\tcat_c$ is a support datum. 
The space $X$ is spectral by assumption, so in order to prove that $(X,\sigma|_{\tcat_c})$ is classifying, we have to show that the assignments
\begin{eqnarray*} 
Y &\mapsto& \mathcal C_{Y}=\{A\in\mathcal T_c \mid \sigma(A)\subset Y\} \\ 
\mathcal C &\mapsto & \sigma(\mathcal C)= \bigcup_{A\in \mathcal C} \sigma(A),
\end{eqnarray*}
define mutually inverse bijections between the set of Thomason subsets $Y\subset X$ and the set of radical thick $\otimes$-ideals $\mathcal C\subset \tcat_c$.

First of all, the two maps are well-defined: the set $\sigma(\mathcal C)$ is a Thomason subset by (S7) (for any subcategory $\mathcal C \subset \tcat_c$) and $\mathcal C_Y$ is a radical thick $\otimes$-ideal by Lemma \ref{lemma:suppZ}~(a) (for any subset $Y\subset X$).

Now, given a thick $\otimes$-ideal $\mathcal C$ in $\tcat_c$, we have the equality $\mathcal C=\langle \mathcal C\rangle_{\otimes}= \mathcal C_{\sigma(\mathcal C)}$ by Lemma \ref{lemma:self_dual} and Lemma \ref{lemma:keycat} applied to $\mathcal E:=\mathcal C$.
Conversely, let $Y=\bigcup_iZ_i$ be a union of closed subsets of~$X$, each with quasi-compact complement $X\smallsetminus Z_i$. Clearly $\sigma(\mathcal C_Y)\subset Y$ by definition (indeed for any subset~$Y\subset X$). By axiom (S8) there are compact objects $A_i$ with $\sigma(A_i)=Z_i$. But then $A_i\in \mathcal C_{Z_i}\subset \mathcal C_Y$, and thus $Y=\bigcup_i \sigma(A_i) \subset \sigma(\mathcal C_Y)$. So we have proved that $\sigma(\mathcal C_Y)=Y$, concluding the verification that the functions $Y\mapsto \mathcal C_Y$ and $\mathcal C\mapsto \sigma(\mathcal C)$ are the inverse of each other.
\qed

\subsection{Compact objects and central rings}
\label{subsec:enter}

In Lemma \ref{lemma:hml_supp} we had ignored conditions (S7) and~(S8). In this section we explore them for the situation when $(X,\sigma)$ can be defined \emph{on compact objects} by functors of the form $\Hom^*_{\tcat}(C,?)_{\mathfrak p}$, where we localize the $\mathrm R_{\tcat}$-module (resp.\ the graded $\mathrm R_{\tcat}^*$-module) $\Hom^*_{\tcat}(C,?)$ with respect to prime ideals $\mathfrak p\in \spec(\mathrm R_{\tcat})$ (resp.\ homogeneous prime ideals $\mathfrak p\in \specgr(\mathrm R^*_{\tcat})$). At a crucial point, we must require that the (graded) central ring is noetherian. Just to be safe, let us explain what we mean precisely by ``localization at a homogeneous prime''.
 
\begin{construction} \label{rem:gr_loc}
Let $M$ be a graded module over a graded commutative ring~$R$. Let $S\subset R$  be a multiplicative system of homogeneous and central elements. Then the localized module
$
S^{-1}M= \{\frac{m}{s} \mid m\in M,s\in S\}
$
is a well-defined graded $S^{-1}R$-module. For a point $\mathfrak p\in \specgr(R)$, we set 
$
M_{\mathfrak p}:= S_{\mathfrak p}^{-1}M,
$
where $S_{\mathfrak p}$ consists of all homogeneous central elements in $R\smallsetminus \mathfrak p$. We write
$
\Supp_R(M)
$
for the `big' support of a graded $R$-module $M$ defined by
$
\Supp_R(M):= \{\mathfrak p\in \specgr(R)\mid M_{\mathfrak p}\not\simeq 0\}
$. 
\end{construction}

For the rest of this section, let $\tcat$ be a compactly generated $\otimes$-triangulated category. Recall from Remark \ref{rem:actions} that the graded Hom sets $\Hom^*_{\tcat}(A,B)$ are graded modules over the graded central ring~$\mathrm R^*_{\tcat}$. We assume given a graded commutative ring $R$ and a grading preserving homomorphism $\phi: R\to \mathrm R^*_{\tcat}$, and always regard the graded Hom sets of $\tcat$ as graded $R$-modules via $\phi$ and the (left) canonical action of~$\mathrm R^*_{\tcat}$. We shall be ultimately interested in the case when $\phi$ is the identity of $\mathrm R^*_{\tcat}$ or the inclusion $\mathrm R_{\tcat}\hookrightarrow \mathrm R^*_{\tcat}$ of its zero degree part (see Prop.\ \ref{prop:S7S8} below). 

\begin{notation} \label{notation:supps}
For each object $A\in\tcat$, define the following subsets of $\specgr(R)$:
\begin{eqnarray*}
\Supptot(A) &:=& \Supp_R(\End^*_{\tcat}(A)) \\
\Supp_B(A) &:=& \Supp_R (\Hom^*_{\tcat}(B,A)) \; , \quad \textrm{for an object }B\in \tcat \\
\Supp_{\mathcal E}(A) &:=& \bigcup_{B\in \mathcal E} \Supp_R (\Hom^*_{\tcat}(B,A)) \; , \quad \textrm{for a family } \mathcal E\subset \tcat.
\end{eqnarray*}
\end{notation}

\begin{lemma} \label{lemma:supps}
In the above notation, we have:
\begin{itemize}
\item[(a)] $\Supptot = \Supp_{\tcat}$.

\item[(b)] Let $E$ be a unital graded $R$-algebra (e.g.\ $E=\End^*_{\tcat}(A)$ for an $A\in\tcat$). Then $\Supp_R(E)= V(\Ann_R(E))$, where the annihilator $\Ann_R(E)$ is the ideal generated by the homogeneous $r\in R$ such that $rE=0$.




\end{itemize}
\end{lemma}

\proof
(a)
Let $A\in \tcat$ and $\mathfrak p\in \specgr(R)$. We have equivalences: $\mathfrak p\not\in \Supptot(A)$
$\Leftrightarrow$ $\id_A=0$ in $\End^*_{\tcat}(A)_{\mathfrak p}$ 
 $\Leftrightarrow$ $f=\id_Af=0$ in $\Hom^*_{\tcat}(B,A)_{\mathfrak p}$ for all $B\in \tcat$ and all $f\in \Hom^*_{\tcat}(B,A)$
  $\Leftrightarrow$ $\mathfrak p\not\in \Supp_{\tcat}(A)$.

(b) Let $\mathfrak p\in \specgr(R)$. Then
$\mathfrak p\not\in V(\Ann_R(E))$
  $\Leftrightarrow$ $\exists$ homogeneous element $r\in R\smallsetminus \mathfrak p$ with $r1_E=0$
    $\Leftrightarrow$ $\exists$  homogeneous central $r\in R\smallsetminus \mathfrak p$ with $r1_E=0$
 (for ``$\Rightarrow$'' simply take $r^2$, which is central because even-graded)
     $\Leftrightarrow$ $E_{\mathfrak p}\simeq0$ 
$\Leftrightarrow$ $\mathfrak p\not\in \Supp_R(E)$.

\qed

\begin{lemma} \label{lemma:equiv_supps}
Let $\mathcal E\subset \tcat$ be a family of objects containing the $\otimes$-unit~$\unit$ and let $X\subset \specgr(R)$ be a subset of homogeneous primes. Assume that the support $(X,\sigma_{X,\mathcal E})$ on $\tcat_c$ defined by $\sigma_{X,\mathcal E}(A):=\Supp_{\mathcal E}(A)\cap X$ satisfies axiom (S5) in Theorem \ref{thm:abstract}, namely: $\sigma_{X,\mathcal E}(A\otimes B)=\sigma_{X,\mathcal E}(A)\cap \sigma_{X,\mathcal E}(B)$ for all $A,B\in \tcat_c$.
Then
\begin{equation*}
\sigma_{X,\mathcal E}(A)= \Supp_{\mathrm{tot}}(A)\cap X
\end{equation*}
for every compact object $A\in \tcat_c$.
\end{lemma}

In particular, if $(X,\sigma_{X,\mathcal E})$ satisfies (S5) then it does not depend on $\mathcal E$!

\proof 
By Lemma \ref{lemma:supps} (a) we have
$$
\sigma_{X,\mathcal E}(A)\stackrel{\mathrm{Def.}}{=} \Supp_{\mathcal E}( A)\cap X \subset \Supp_{\tcat}(A)\cap X = \Supptot(A)\cap X
$$
for all~$A$.
By our convention, every compact object in $\tcat$ is rigid. 
It follows that
\begin{eqnarray*}
\Supptot(A) \cap X &=&
 \Supp_A(A) \cap X \\
&\stackrel{A \textrm{ rigid}}{=}& 
    \Supp_{\unit}(A^{\vee}\otimes A) \cap X \\
&=& 
 \sigma_{X,\{\unit\}}(A^{\vee}\otimes A) \\
&\subset & 
 \sigma_{X,\mathcal E }(A^{\vee}\otimes A) \\
&\stackrel{\textrm{(S5)}}{=}&
 \sigma_{X,\mathcal E } (A^{\vee})\cap\sigma_{X,\mathcal E } (A) \\ 
 &\subset&
  \sigma_{X,\mathcal E }(A),
\end{eqnarray*}
thus proving the reverse inclusion.
\qed

\begin{prop} \label{prop:S7S8}
Let $\tcat$ be a compactly generated $\otimes$-triangulated category. Let $R$ be either the graded central ring $\mathrm R^*_{\tcat}$ or its subring $\mathrm R_{\tcat}$, and assume that it is (graded) noetherian. Let $(X,\sigma_{X}:=\sigma_{X,\{\unit\}})$ be the support on $\tcat_c$ we defined in Lemma \ref{lemma:equiv_supps}, for some subset $X\subset \specgr(R)$, and again assume that $(X,\sigma_{X})$ satisfies (S5) on $\tcat_c$. Then
\begin{itemize}
\item[(a)] The support $(X,\sigma_X)$ satisfies axiom (S7) in Theorem \ref{thm:abstract}, namely: For every $A\in \tcat_c$ the subset $\sigma_{X}(A)$ is closed in $X$ and its complement $X\smallsetminus \sigma_X(A)$ is quasi-compact. 

\item[(b)] The support $(X,\sigma_X)$ satisfies axiom (S8) in Theorem \ref{thm:abstract}: For every closed subset $Z\subset X$ there exists a compact object $A\in \tcat_c$ with $\sigma_{X}(A)=Z$. 

\end{itemize}
\end{prop}

\proof
(a) By Lemma \ref{lemma:equiv_supps} and Lemma \ref{lemma:supps} (b), for each $A\in \tcat_c$ we have equalities
$$\sigma_X(A)=\Supp_{\mathrm{tot}}(A)\cap X=V(\Ann_R(\End^*_{\tcat}(A)))\cap X.$$
This is by definition a closed subset of~$X$. Since we assumed $R$ noetherian, it follows easily that \emph{every} open subset of $\specgr(R)$ is quasi-compact.

(b)
Every closed subset of $X$ has the form $Z=X\cap V(I)$ for some homogeneous ideal $I<R$. Since $R$ is noetherian, $I$ is generated by finitely many homogeneous elements, say $I=\langle r_1,\ldots, r_n \rangle$.
Let $C_i$ be the cone of $r_i:\unit\to T^{m_i}\unit $. It is rigid and compact, and moreover we claim that 
$\Supp_{\unit}(C_i)=V(\langle r_i\rangle )$. 
Indeed, by applying $\Hom_{\tcat}^*(\unit,\,?)_{\mathfrak p}$ to the distinguished triangle $\unit\stackrel{r_i}{\to} T^{m_i}\unit\to C_i\to T\unit$, we obtain an exact sequence
\begin{equation*}
\xymatrix@1{
\Hom_{\tcat}^*(\unit,\unit)_{\mathfrak p} \ar[r]^-{ r_i\cdot } &  \Hom_{\tcat}^{*+m_i}(\unit,\unit)_{\mathfrak p} \ar[r] & \Hom_{\tcat}^*(\unit,C_i)_{\mathfrak p} \ar[r] & \Hom_{\tcat}^{*+1}(\unit,\unit)_{\mathfrak p}
}
\end{equation*}
of graded $R$-modules. 
Note that the first morphism is multiplication by $r_i$ (see \ref{rem:actions}). It is invertible if \emph{and only if} $r_i$ is invertible in $R_{\mathfrak p}$, because we assumed that $R=\mathrm R^*_{\tcat}$ or $R=\mathrm R_{\tcat}$.
Hence $r_i\in R^{\times}_{\mathfrak p}$ $\Leftrightarrow$ $\Hom_{\tcat}^*(\unit,C_i)_{\mathfrak p}\simeq 0$ $\Leftrightarrow$ $\mathfrak p\not\in \Supp_{\unit}(C_i)$, as claimed.
Now it suffices to set $A:=C_1\otimes\cdots\otimes C_n$ (which is again a rigid and compact object by Conv.\ \ref{conv:cptly_gen}~(b)), because then
\begin{eqnarray*}
\sigma_{X}(A) 
&\stackrel{\textrm{(S5)}}{=}&
 \sigma_{X}(C_1)\cap \cdots \cap \sigma_{X}(C_n) \\
&=&
 X\cap \Supp_{\unit}(C_1)\cap\cdots \cap \Supp_{\unit}(C_n) \\
&= &
 X\cap V(\langle r_1\rangle)\cap\cdots \cap V(\langle r_n\rangle)   \\
&=&
 X\cap V(I)\quad=\quad Z\;,
\end{eqnarray*}
as desired. 
\qed

\subsection{Comparison with the support of Benson-Iyengar-Krause}
\label{subsec:comparison}

As an application of the last two sections, we provide sufficient conditions for the support defined by Benson, Iyengar and Krause in \cite{bik} to coincide with Balmer's support on compact objects, in the situation where both supports are defined.

Let $\tcat$ be a tensor triangulated category which is a genuine compactly generated category, such that the tensor is exact and preserves small coproducts in both variables, and where compact and rigid objects coincide (thus in particular $\tcat$ satisfies the hypotheses in Convention \ref{conv:cptly_gen}). Let $R$ be either $\mathrm R^*_{\tcat}= \End^*_{\tcat}(\unit)$ or $\mathrm R_{\tcat}=\End_{\tcat}(\unit)$, and assume that it is a (graded) noetherian ring. In such a situation, the support $\supp_R^{\mathrm{BIK}}: \obj(\tcat)\to 2^{\specgr(R)}$ defined in \cite{bik} can be given by the formula 
\begin{equation} \label{supp_bik}
\supp_R^{\mathrm{BIK}}(A)=\{\mathfrak p \mid \mathit{\Gamma}_{\mathfrak p}(\unit)\otimes A\not\simeq 0\}
\quad \subset \quad \specgr(R) 
\end{equation}
for every $A\in \tcat$, where $\mathit \Gamma_{\mathfrak p}(\unit)$ is a certain non-trivial object depending on~$\mathfrak p$ (see \emph{loc.\ cit.}, especially $\S5$ and Cor.\ 8.3). In this setting, $\supp_R^{\mathrm{BIK}}$ also recovers the support for noetherian stable homotopy categories considered in \cite[$\S$6]{hps}. 

Here is our comparison result:

\begin{thm}\label{prop:bik_class} 
 Keep the notation of the last paragraph. Let further $X\subset \specgr(R)$ be a spectral subset, and write $\sigma(A):=X\cap \supp_R^{\mathrm{BIK}}(A)$ for the restricted support. Assume the following three hypotheses:
\begin{itemize}
\item[(1)]  For every compact $A\in \tcat_c$, we have $\sigma(A)=X\cap V(\Ann_R (\End^*_{\tcat}(A)))$.

\item[(2)] The support $(X,\sigma)$ detects objects of $\tcat$:  $\sigma(A)=\emptyset \Rightarrow A\simeq0$.

\item[(3)] The support $(X,\sigma)$ satisfies the `partial Tensor Product theorem':
$$
\sigma(A\otimes B)= \sigma(A)\cap \sigma(B)
$$
whenever $A\in \tcat_c$ is compact and $B\in \tcat$ arbitrary. 

\end{itemize}
Then there is a unique isomorphism $(X,\sigma)\simeq (\spc(\tcat_c), \supp)$ of support data on $\tcat_c$ between the restricted Benson-Iyengar-Krause support and the Balmer support.
\end{thm}

\begin{remark}
Note that hypothesis (1) is not so restrictive as it may seem. Indeed, by \cite[Thm.\ 5.5]{bik} it must hold for every $A\in \tcat_c$ for which $\End^*_{\tcat}(A)$ is finitely generated over~$R$. Also, (2) holds for the choice $X:= \specgr(R)$ by \cite[Thm.\ 5.2]{bik}. Thus, our theorem says roughly that, if we can `adjust' the Benson-Iyengar-Krause support by restricting it to a subset $X$ in such a way that it satisfies  the partial Tensor Product theorem and it still detects objects, then it must be the universal support datum on~$\tcat_c$. 
\end{remark}

\proof
It suffices to show that $(X,\sigma)$ satisfies axioms (S0)-(S9) in Theorem \ref{thm:abstract}. Note that 
(S0)-(S4) and (S6) are immediate from  \eqref{supp_bik}, and (S5), resp.\ (S9), are simply assumed in hypothesis (3), resp.~(2). 
We are left with the verification of (S7) and (S8). By hypothesis (1), the restriction of $(X,\sigma)$ on compact objects coincides with the support $(X,\sigma_X)=(X,\sigma_{X,\mathcal E})$ of the previous section $\S$\ref{subsec:enter}. Hence, since $R$ is noetherian, $(X,\sigma)$ satisfies (S7) and (S8) by virtue of Proposition \ref{prop:S7S8}.
\qed

\section{The spectrum and the Baum-Connes conjecture}
\label{sec:BC}

As in the Introduction, let $G$ be a second countable locally compact Hausdorff group, and let $\KK^G$ be the $G$-equivariant Kasparov category of separable $G$-$\Cstar$-algebras (see \cite{meyernest-bc} \cite{mayer-cat}). It is a tensor triangulated category as in Definition \ref{defi:tensor_tr}, with arbitrary countable coproducts (\cite[App.\ A]{meyernest-bc} \cite[App.\ A]{thesis}). The tensor structure $\otimes$ is induced by the minimal tensor product of $\Cstar$-algebras with the diagonal $G$-action, and the unit object $\unit$ is the field of complex numbers $\C$ with the trivial $G$-action. Of the rich functoriality of $\KK^G$, we mention the \emph{restriction} tensor triangulated functor
$
\Res{G}{H}: \KK^G \to \KK^H
$
and the \emph{induction} triangulated functor
$
\Ind{H}{G}: \KK^H \to \KK^G
$
for $H$ a closed subgroup of~$G$.
They are related by a `Frobenius' natural isomorphism
\begin{equation} \label{frobenius}
\Ind{H}{G}(A\otimes \Res{G}{H}(B))\simeq \Ind{H}{G}(A) \otimes B.
\end{equation}
Roughly speaking, the Baum-Connes Conjecture proposes a computation for the $K$-theory of the \emph{reduced crossed product} $G\ltimes \,?:\KK^G\to \KK$.
We recall now the conceptual formulation of the conjecture, and its generalizations, due to Meyer and Nest \cite{meyernest-bc}.

\begin{defi}
Consider the two full subcategories of $\KK^G$
$$\CI^G:= \bigcup_{H\leq G \textrm{ compact}} \Image(\Ind{H}{G}) 
\quad \textrm{ and } \quad
\CC^G:= \bigcap_{H\leq G \textrm{ compact}}\Ker(\Res{G}{H})$$
(for ``compactly induced'' and ``compactly contractible'', respectively). We consider the localizing hull $\langle \CI^G\rangle_{\loc}\subset \KK^G$. Note that both $\langle \CI^G\rangle_{\loc}$ and $\CC^G$ are  localizing subcategories. Both are also $\otimes$-ideals: $\CC^G$ because each $\Res{G}{H}$ is a $\otimes$-triangulated functor and $\langle\CI^G\rangle_{\loc}$ because of the Frobenius formula \eqref{frobenius}.
\end{defi}

\begin{thm}[{\cite[Thm.\ 4.7]{meyernest-bc}}]
The localizing tensor ideals $\langle \CI^G\rangle_{\loc}$ and $\CC^G$ are complementary in~$\KK^G$ (see Def.~\ref{defi:compl_pair}).
\qed
\end{thm}

By Remark \ref{remark:tensor_gluing}, the gluing triangle for this complementary pair at any object $A\in\KK^G$, that we shall denote by $P^G(A)\stackrel{D^G(A)}{\to} A \to N^G(A) \to TP^G(A)$,  is obtained by tensoring $A$ with the gluing triangle
\begin{equation*}
\xymatrix@1{
P^G(\unit) \ar[r]^-{D^G(\unit)} & \unit \ar[r] & N^G(\unit) \ar[r] & TP^G(\unit)
}
\end{equation*}
 for the tensor unit. The approximation $D^G=D^G(\unit): P^G(\unit)\to \unit$ is called the \emph{Dirac morphism for $G$}. Note that, by the general properties of Bousfield localization (Prop.\ \ref{prop:compl}), the objects $P^G(\unit)$ and $N^G(\unit)$ are $\otimes$-idempotent: 
\begin{equation} \label{tensidem}
P^G(\unit)\otimes P^G(\unit)\simeq P^G(\unit) \quad, \quad N^G(\unit)\otimes N^G(\unit)\simeq N^G(\unit).
\end{equation}

\begin{defi}\label{defi:BC}
Let $A\in \KK^G$, and let $F:\KK^G\to \mathcal C$ be any functor defined on the equivariant Kasparov category. One says that $G$ \emph{satisfies the Baum-Connes conjecture for $F$ with coefficients $A$} if the homomorphism
\begin{equation} \label{assembly}
F( D^G(A))\;:\; F( P^G(A))\longrightarrow F(A)
\end{equation}
is an isomorphism in~$\mathcal C$. 
\end{defi}

The main result of \cite{meyernest-bc} is a proof that, if $F=K_*(G\ltimes\, ?): \KK^G\to \Ab$ is the $K$-theory of the reduced crossed product, then the homomorphism \eqref{assembly} is naturally isomorphic to the so-called assembly map for the group $G$ with coefficients~$A$, implying that for this choice of $F$ the above formulation of the Baum-Connes conjecture is equivalent to the original formulation with coefficients (see \cite{bch}).

The above formulation for general functors $F$ on $\KK^G$ is then a natural generalization.  
Note that, if the Dirac morphism $D^G$ is itself an isomorphism in $\KK^G$, then  $G$ satisfies the conjecture for all functors $F$ and all coefficients $A\in \KK^G$. Note also that $D^G$ is an isomorphism if and only if $N^G(\unit)\simeq 0$, if and only if the inclusion $\langle \CI^G\rangle_{\loc}\hookrightarrow \KK^G$ is an equivalence.

In \cite{higson-kasparov}, Higson and Kasparov proved that the Dirac morphism is invertible, and therefore that the conjecture holds for every functor and all coefficients, for groups $G$ having the \emph{Haagerup approximation property} ($=$ \emph{a-T-menable} groups). These are groups admitting a proper and isometric action on Hilbert space, in a suitable sense. They form a rather large class containing all amenable groups.

 We contribute the following intriguing observation, which serves as a motivation for pursuing the (tensor triangular) geometric study of triangulated categories arising in connection with Kasparov theory.
\begin{thm} \label{thm:covering-bc}
Assume that  the spectrum of $\KK^G$ is covered by the spectra of $\KK^H$ as $H$ runs through the compact subgroups of~$G$:
\begin{equation} \label{covhyp}
\spc\big( \KK^G \big) = \bigcup_{H\leq G \textrm{ compact}} \Spc(\Res{G}{H}) \Big( \spc\big(\KK^H \big)\Big).
\end{equation}
Then the Dirac morphism $D^G:P^G(\unit)\to \unit$ is an isomorphism.
\end{thm}

\proof By a basic result of tensor triangular geometry (see \cite[Cor.\ 2.4]{balmer_prime}), an object $A\in\KK^G$ belongs in each prime $\otimes$-ideal $\mathcal P\in\spc(\KK^G)$ if and only if it is $\otimes$\hyph nilpotent, i.e., if and only if $A^{\otimes n}\simeq 0$ for some $n\geq1$. Thus if the covering hypothesis \eqref{covhyp} holds, we have
\begin{eqnarray*}
A\,\textrm{ is }\otimes\textrm{-nilpotent }&\Leftrightarrow & A\in \mathcal P \quad \forall \mathcal P\in\spc(\KK^G) \\
&\Leftrightarrow & A \in(\Res{G}{H})^{-1}\mathcal Q \quad \forall \mathcal Q\in \spc(\KK^H),\forall \,H \\
&\Leftrightarrow & \Res{G}{H}(A) \in\mathcal Q \quad \forall \mathcal Q\in \spc(\KK^H),\forall \,H
\end{eqnarray*}
where $H$ ranges among all compact subgroups of $G$. Now specialize the above to $A:=N^G(\unit)$. Clearly $N^G(\unit)$ satisfies the latter condition, because by construction $N^G(\unit)\in \CC^G=\bigcap_{H}\Ker(\Res{G}{H})$. Thus $N^G(\unit)$ is a $\otimes$-nilpotent object. But $N^G(\unit)$ is also $\otimes$-\emph{idem}potent \eqref{tensidem}, and therefore $N^G(\unit)\simeq 0$, implying the claim.
\qed

\section{Some homological algebra for $KK$-theory}

 We recall a few definitions and results of relative homological algebra in triangulated categories (\cite{christensen} \cite{beligiannis} \cite{meyernest-hom}); our reference is \cite{meyernest-hom}. Here $\tcat$ will always denote a triangulated category admitting (at least) all countable coproducts.

\begin{defi} \label{defi:stable}
 A \emph{stable abelian category} is an abelian category $\mathcal A=(\mathcal A,T)$ equip\-ped with a self-equivalence $T:\mathcal A\stackrel{\sim}{\to} \mathcal A$. A \emph{stable homological functor} $H=(H,\delta)$ on $\tcat$ is an additive functor $H:\tcat \to \mathcal A$ to some stable abelian category $\mathcal A$ together with an isomorphim $\delta:HT\stackrel{\sim}{\to}TH$, and such that for every distinguished triangle $A \stackrel{u}{\to}  B \stackrel{v}{\to}  C \stackrel{w}{\to}  TA$ of $\tcat$ the sequence $HA \stackrel{Hu}{\to}  HB \stackrel{Hv}{\to} HC \stackrel{\delta Hw}{\longrightarrow}  THA$ is exact in~$\mathcal A$. 
\end{defi}

\begin{example} \label{ex:hlg}
If $H:\tcat\to \mathcal A$ is a homological functor in the usual sense (i.e., an additive functor to some abelian category $\mathcal A$ such that if $A\to B\to C\to TA$ is distinguished in $\tcat$ then $HA\to HB\to HC$ is exact), we may construct a stable homological functor $H_{*}:\tcat\to \mathcal A^{\Z}$ as follows. Let $\mathcal A^{\Z}$ be the category of $\Z$-graded objects $M_*=(M_n)_{n\in\Z}$ in $\mathcal A$ (with degree-zero morphisms); with the shift $TM_*:=(M_{n-1})_n$ it is a stable abelian category. Then $H_*(A):=(HT^{-n}A)_n$ defines a stable homological functor (with $\delta=\id$).
Note that, if the translation $T$ of $\tcat$ is $n$-periodic for some $n\geq 1$, by which we mean that there is an isomorphism $T^n\simeq \id_{\tcat}$, then we may equally consider $H_*$ as a functor to the stable abelian category $\mathcal A^{\Z/n}$ of $\Z/n$-graded objects of~$\mathcal A$.
\end{example}
 
 \begin{defi} \label{defi:hlg_alg} 
 A \emph{homological ideal} $\mathcal I$ in $\tcat$ is a subfunctor $\mathcal I\subset \Hom_{\tcat}(\textrm{?`},?)$ of the form $\mathcal I=\ker(H)$ for some stable homological functor~$H:\tcat\to \mathcal A$. 
 For convenience, we define a \emph{homological pair} $(\tcat,\mathcal I)$ to consist of a triangulated category $\tcat$ with countable coproducts together with a homological ideal $\mathcal I$ in $\tcat$ which is closed under the formation of countable coproducts of morphisms. If $\mathcal I=\ker(H)$, the last condition is satisfied whenever $H$ commutes with countable coproducts. 
 
 Let $(\tcat,\mathcal I)$ be a homological pair. A (stable) homological functor $H:\tcat\to \mathcal A$ is \emph{$\mathcal I$-exact} if $H(f)=0$ for all $f\in\mathcal I$.
 An object $P\in \tcat$ is \emph{$\mathcal I$-projective} if $\Hom(P,?):\tcat\to \Ab$ is $\mathcal I$-exact. 
 An object $N\in \tcat$ is \emph{$\mathcal I$-contractible} if $\id_N\in \mathcal I$. The category $\tcat$ \emph{has enough $\mathcal I$-projectives} if, for every $A\in \tcat$, there exists a distinguished triangle $B\to P\to A \to TB$ where $P$ is $\mathcal I$-projective  and $(A\to TB)\in \mathcal I$. 
\end{defi}

\begin{remark} \label{rem:univ_exact}
It can be shown that for every pair $(\tcat,\mathcal I)$ there exists a universal $\mathcal I$-exact stable homological functor $h_{\mathcal I}:\tcat\to \mathcal A(\tcat,\mathcal I)$ (where $\mathcal A(\tcat, \mathcal I)$ has small hom sets) --  at least if $\tcat$ has enough $\mathcal I$-projectives, which is the case in all our examples. See \cite[$\S$3.7]{meyernest-hom} for details. 
With this assumption, it is proved in \emph{loc.\ cit.} that $h_{\mathcal I}$ restricts to an equivalence between the full subcategory $\mathcal P_{\mathcal I}$ of $\mathcal I$-projective objects in $\tcat$ and the full subcategory of projectives in the stable abelian category~$\mathcal A(\tcat,\mathcal I)$.
\end{remark}

\begin{thm}[{\cite[Thm.\ 3.16]{meyer-hom}}] \label{thm:hl_compl_pair}
Let $(\tcat,\mathcal I)$ be a homological pair, and assume that $\mathcal T$ has enough $\mathcal I$-projectives. Then the pair of subcategories $(\langle \mathcal P_{\mathcal I}  \rangle_{\loc} , \mathcal  N_{\mathcal I})$ is complementary in~$\tcat$, where $\mathcal P_{\mathcal I}$ denotes the full subcategory of $\mathcal I$-projective objects in $\tcat$ and $\mathcal N_{\mathcal I}$ that of $\mathcal I$-contractible ones.
\end{thm}

Fix a homological pair $(\tcat,\mathcal I)$. Given additive functors $F:\tcat\to \mathcal C$ and $G:\tcat^{\op}\to \mathcal D$ to some abelian categories $\mathcal C,\mathcal D$, if there are enough $\mathcal I$-projective objects one may use $\mathcal I$-projective resolutions to define, in the usual way, both the \emph{left derived functors} $\mathsf L_n^{\mathcal I}F:\tcat\to \mathcal C$ and the \emph{right derived functors} $\mathsf R_{\mathcal I}^nG:\tcat^{\op}\to \mathcal D$ (\emph{relative to $\mathcal I$}), for $n\geq 0$. These can sometimes be identified with more familiar derived functors in the context of abelian categories by means of the universal exact functor $h_{\mathcal I}:\tcat \to \mathcal A(\tcat,\mathcal I)$ (see e.g.\ Prop.\ \ref{prop:identif} below).  The notation $\Ext_{\tcat,\mathcal I}^n(A,B)$ stands for $\mathsf R_{\mathcal I}^nG(A)$ in the case of the functor $G=\Hom_{\tcat}(\textrm{?`},B):\tcat^{\op}\to \Ab$.

We will make use of some instances of the following result:

\begin{thm}
\label{thm:hereditary}
Let $(\mathcal T,\mathcal I)$ be a homological pair. Let $A\in \langle \mathcal P_{\mathcal I}\rangle_{\loc}$, and assume that $A$ has an $\mathcal I$-projective resolution of length one. Then
\begin{itemize}
\item[(a)] 
For every homological functor $F:\tcat\to \mathcal A$ there is a natural exact sequence
\begin{equation*}
\xymatrix@1{
0 \ar[r] & \mathsf L_0^{\mathcal I}F(A) \ar[r] & F(A) \ar[r] & \mathsf L_1^{\mathcal I}F(TA)\ar[r] &0.
}
\end{equation*}
\item[(b)]
For every homological functor $G:\tcat^{\op}\to \mathcal A$ 
there is a natural exact sequence
\begin{equation*}
\xymatrix@1{
0 \ar[r] & \mathsf R^1_{\mathcal I}G(TA) \ar[r] & G(A) \ar[r] & \mathsf R^0_{\mathcal I}G(A)\ar[r] &0.
}
\end{equation*}
\item[(c)]
Choosing $G=\Hom_{\tcat}(\textrm{\emph{?`}},B)$ in (b), for any object $B\in \tcat$, we get
\begin{equation*}
\xymatrix@1{
0 \ar[r] & \Ext^1_{\mathcal T,\mathcal I}(TA,B) \ar[r] & \Hom_{\tcat}(A,B) \ar[r] & \Ext^0_{\mathcal T,\mathcal I}(A,B)\ar[r] &0.
}
\end{equation*}
\end{itemize}
\end{thm}

\proof
This is {\cite[Thm.\ 66]{meyernest-hom}}. Note that our assumption $A\in \langle \mathcal P_{\mathcal I}\rangle_{\loc}$ coincides with that in \emph{loc.\ cit.}, namely $A\in {}^{\perp}\mathcal N_{\mathcal I}$, because of Theorem \ref{thm:hl_compl_pair}.
\qed

\begin{remark} \label{rem:splitting}
In the situation of Theorem \ref{thm:hereditary}, assume that there exists a decomposition $A\cong A_0\oplus A_1$ such that $\mathsf L_i^{\mathcal I}F(A_j)=0$ (resp.\ $\mathsf R^i_{\mathcal I}G(A_j)=0$) for $\{i,j\}=\{0,1\}$. Then we see from its naturality and additivity that the sequence in (a) (resp.\ in (b) and (c)) has a splitting, determined by the isomorphim $A\simeq A_0\oplus A_1$.
\end{remark}

\subsection{The categories $\tcat^G$ and $\mathcal K^G$}
\label{subsec:cats}

Consider the equivariant Kasparov category $\KK^G$ for a compact group~$G$. We recall that the $R(G)$-modules $\Hom_{\KK^G}(T^i\unit, A)=\KK^G(T^i\unit, A)$ identify naturally with topological $G$-equivariant $K$-theory $K^G_i(A)$ (\cite[$\S$2]{phillips-free}, \cite[$\S$11]{blackadar}). 
By the Green-Julg theorem (\cite[Thm.\ 11.7.1]{blackadar}), there is an isomorphism $K^G_i\simeq K_i(G\ltimes ?)$. Since ordinary $K$-theory $K_*$ of separable $\Cstar$-algebras yields countable abelian groups and commutes with countable coproducts in $\KK^G$, and since $G\ltimes ?$ commutes with coproducts and preserves separability, we conclude that the $\otimes$-unit $\unit=\C$ is a compact${}_{\aleph_1}$ object of  $\KK^G$ (Def.\ \ref{defi:cpt}). Hence the category
$\tcat^G:= \langle \unit \rangle_{\loc} \subset \KK^G$ is compactly${}_{\aleph_1}$ generated. Moreover, since  it is monogenic -- in the sense of being generated by the translations of the $\otimes$-unit -- its compact and rigid objects coincide, and form a thick  $\otimes$-triangulated subcategory $\mathcal K^G:=\tcat^G_c=\langle \unit \rangle$, which is also the smallest thick subcategory of $\KK^G$ containing the tensor unit.  In particular $\tcat^G$ is a compactly generated $\otimes$-triangulated category as in Convention \ref{conv:cptly_gen}. 

 As in $\KK^G$, we have Bott periodicity: $T^2\simeq \id_{\tcat^G}$. Hence all homological functors $H:\tcat^G\to \mathcal A$ give rise to stable homological functors $H_*$ to the category of $\Z/2$-graded objects~$\mathcal A^{\Z/2}$ (see Example \ref{ex:hlg}).

The relevance of $\tcat^G$ to $K$-theory is explained by the following result.

\begin{thm} \label{thm:compl_Ktheory}
Let $G$ be a compact group. The pair of localizing subcategories $(\tcat^G, \Ker (K^G_*))$ of $\KK^G$ is complementary. In particular, there exists a triangulated functor $L:\KK^G\to \tcat^G$ and a natural map $L(A)\to A$ inducing an isomorphism $K^G_*(LA)\simeq K^G_*(A)$ for all $A\in \KK^G$.
\end{thm}

\proof 
Meyer and Nest prove (\cite[Thm.\ 72]{meyernest-hom}) that $K^G_*= K_*\circ (G\ltimes ?)$, as a functor from $\KK^G$ to $\Z/2$-graded countable $R(G)$-modules, is the universal $\ker(K^G_*)$-exact functor and that, as a consequence, it induces an equivalence between the category $\mathcal P_{\ker(K^G_*)}$ of $\ker(K_*^G)$-projective objects in $\KK^G$ and that of projective graded $R(G)$-modules (cf.\ Remark \ref{rem:univ_exact}). Since every projective module is a direct summand of a coproduct of copies of $R(G)=K^G_*(\unit)$ and of its shift $R(G)(1)=K^G_*(T\unit)$, it follows that $\langle \mathcal P_{\ker(K^G_*)}\rangle_{\loc} = \langle \unit \rangle_{\loc} \subset \KK^G $, and therefore the claim is just  Theorem \ref{thm:hl_compl_pair} applied to the homological pair $(\KK^G, \ker(K^G_*))$.
\qed

We shall make use of quite similar arguments in the following section.

 In the rest of this article we shall begin the study of these categories from a geometric point of view, concentrating on the easier case of a finite group~$G$.


\subsection{Central localization of equivariant $KK$-theory}

\label{sec:locKth}

Let $G$ be a compact group, and let $\mathfrak p\in\spec(R(G))$.
We wish to apply the abstract results of $\S$\ref{subsec:central} to the monogenic compactly generated tensor triangulated category $\tcat=\tcat^G$ and the  multiplicative system $S=R(G)\smallsetminus \mathfrak p$. Thus we consider the thick $\otimes$-ideal of compact objects
\begin{equation*}
\mathcal J^G_{\mathfrak p}:=\langle \cone(s)\mid s\in R(G)\smallsetminus \mathfrak p\rangle_{\otimes}\subset \tcat^G_c
\end{equation*}
and the localizing $\otimes$-ideal
$ 
\mathcal L^G_{\mathfrak p}:=\langle \mathcal J^G_{\mathfrak p} \rangle_{\loc} \subset \tcat^G
$ 
that it generates. We denote its right orthogonal category of $\mathfrak p$-local objects by
\begin{equation}
\tcat^G_{\mathfrak p}:=(\mathcal L^G_{\mathfrak p})^{\perp} \simeq \tcat^G/\mathcal L^G_{\mathfrak p}.
\end{equation}
Now Theorem \ref{thm:lochom} specializes to the following result, which says that $\tcat^G_{\mathfrak p}$ is a well-behaved notion of localization of $\tcat^G$ at~$\mathfrak p$. Note that similar results are true with, instead of $\tcat^G$, any other localizing $\otimes$-subcategory of $\KK^G$ generated by compact and rigid objects, and also, obviously, for multiplicative subsets which do not necessarily come from prime ideals.

\begin{thm} \label{thm:locKth}
The pair 
$
(\mathcal L^G_{\mathfrak p} , \tcat^G_{\mathfrak p} )
$
is a complementary pair of localizing $\otimes$-ideals of~$\tcat^G$. In particular, the gluing triangle for an object $A\in\tcat^G$ is obtained by tensoring $A$ with the gluing triangle for the $\otimes$-unit, which we denote by
\begin{equation} \label{gluing_tri_unit}
\xymatrix{
{}_{\mathfrak p}\unit \ar[r]^-{\varepsilon} & \unit \ar[r]^-{\eta} & \unit_{\mathfrak p} \ar[r] & T({}_{\mathfrak p}\unit) .
}
\end{equation}
Moreover, the following hold true:
\begin{itemize}
\item[(a)] $\mathcal L^G_{\mathfrak p}={}_{\mathfrak p}\unit\otimes \tcat^G$ and  $\mathcal T^G_{\mathfrak p}= \unit_{\mathfrak p} \otimes \tcat^G$.
\item[(b)] The maps $\varepsilon$ and $\eta$ induce isomorphisms ${}_{\mathfrak p}\unit\simeq {}_{\mathfrak p}\unit\otimes {}_{\mathfrak p}\unit $ and $\unit_{\mathfrak p}\simeq \unit_{\mathfrak p} \otimes {\unit}_{\mathfrak p}$.
\item[(c)] The category $\tcat^G_{\mathfrak p}$ is a monogenic compactly generated $\otimes$-tri\-an\-gu\-lat\-ed category with tensor unit $\unit_{\mathfrak p}$.
\item[(d)] Its tensor triangulated subcategory of compact and rigid objects is $(\tcat^G_{\mathfrak p})_c=\langle \unit_{\mathfrak p}\otimes \tcat^G_c\rangle\subset \tcat^G_{\mathfrak p}$. 
\item[(e)] The functor
$
 \unit_{\mathfrak p}\otimes\, ?:\tcat^G\to \mathcal \tcat^G_{\mathfrak p}
$
is an $R(G)$-linear $\otimes$-triangulated functor commuting with coproducts.

\item[(f)] The central ring $\mathsf R_{\tcat^G_{\mathfrak p}}=\End(\unit_{\mathfrak p})$ of $\tcat^G_{\mathfrak p}$ is $R(G)_{\mathfrak p}$, and  $K^G_0(\eta:\unit \to \unit_{\mathfrak p})$ is the localization homomorphism $R(G)\to R(G)_{\mathfrak p}$.

\item[(g)] $A$ is $\mathfrak p$-local (i.e.,\ $A\in \tcat^G_{\mathfrak p}$) $\Leftrightarrow$ $s\cdot \id_A$ is invertible for every $s\in R(G)\smallsetminus {\mathfrak p}$.
\item[(h)] If $A\in \tcat^G_c$, then $\eta:B\to \unit_{\mathfrak p}\otimes B$ induces a canonical isomorphism
\begin{equation*}
\KK^G(A,B)_{\mathfrak p}\simeq \KK^G(A,\unit_{\mathfrak p}\otimes B)
\end{equation*}
for every $B\in \tcat^G$. In particular $K^G_*(B)_{\mathfrak p}\simeq K^G_*(\unit_{\mathfrak p}\otimes B) $ (set $A=T^*\unit$).
\end{itemize}
\end{thm}

\begin{cor} \label{cor:local_approx}
For $G$ a compact group and $\mathfrak p\in \spec(R(G))$, there exist a triangulated functor $L_{\mathfrak p}:\KK^G\to \tcat^G_{\mathfrak p}$ on the equivariant Kasparov category and natural maps $L_{\mathfrak p}(A)\leftarrow L(A)\to A$ in $\KK^G$, inducing an isomorphism $K^G_*(L_{\mathfrak p} A)\simeq K^G_*(A)_{\mathfrak p}$.
\end{cor}

\proof
By Theorem \ref{thm:compl_Ktheory}, there exists in $\KK^G$ a natural map $LA\to A$ with $LA\in \tcat^G$ and $K^G_*(LA \to A)$ invertible. Set $LA \to L_{\mathfrak p}A$ to be $\eta: LA \to \unit_{\mathfrak q}\otimes LA$ as in Theorem \ref{thm:locKth}. The fraction $L_{\mathfrak p}(A)\leftarrow L(A)\to A$ in $\KK^G$ has the required property.
\qed

For later use, we record the behaviour of central localization under restriction.

\begin{lemma} \label{lemma:central_restriction}
Let $H$ be a closed subgroup of the compact group $G$. Moreover, let $\mathfrak q $ be a prime ideal in $R(H)$ and let $\mathfrak p:= (\Res{G}{H})^{-1}(\mathfrak q)\in \spec(R(G))$. Let ${}_{\mathfrak p}\unit\to \unit \to \unit_{\mathfrak p}\to T({}_{\mathfrak p}\unit)$ be the  gluing triangle in $\tcat^G$ for~$\mathfrak p$ and let
${}_{\mathfrak q}\unit\to \unit \to \unit_{\mathfrak q}\to T({}_{\mathfrak q}\unit)$ be the one in $\tcat^H$ for~$\mathfrak q$. Then
\begin{equation*}
\Res{G}{H}({}_{\mathfrak p}\unit)\otimes {}_{\mathfrak q}\unit \simeq \Res{G}{H}({}_{\mathfrak p}\unit)
\quad \textrm{ and } \quad 
\unit_{\mathfrak q}\otimes \Res{G}{H}(\unit_{\mathfrak p}) \simeq \unit_{\mathfrak q}.
\end{equation*}
\end{lemma}

\proof
Note that $S:= \Res{G}{H}(R(G)\smallsetminus \mathfrak p)$ is a multiplicative system in $R(H)$, so there is an associated central localization of $\tcat^H$ with complementary pair $(\mathcal L_{S}^H,\tcat^H_{S})$ and gluing triangle
${}_S\unit\to \unit \to \unit_S \to T({}_S\unit)$. We claim that this triangle is isomorphic to the restriction of ${}_{\mathfrak p}\unit\to \unit \to \unit_{\mathfrak p}\to T({}_{\mathfrak p}\unit)$. By the uniqueness of gluing triangles and since $\Res{G}{H}(\unit)=\unit$, it suffices to show that $\Res{G}{H}(\mathcal L^G_{\mathfrak p})\subset \mathcal L^H_S$ and $\Res{G}{H}(\mathcal T^G_{\mathfrak p})\subset \tcat^H_S$.
The first inclusion holds because $\Res{G}{H}$ is a coproduct preserving $\otimes$-triangulated functor and because $\Res{G}{H}(\cone(s))\simeq \cone(\Res{G}{H}(s))\in \mathcal L^H_S$ for all $s\in R(G)\smallsetminus \mathfrak p$. The second inclusion holds by the characterization in Theorem \ref{thm:lochom} (g) of the objects of $\tcat^H_S$. Finally, the inclusion $S\subset R(H)\smallsetminus \mathfrak q$ implies $\mathcal L^H_S\subset \mathcal L^H_{\mathfrak q}$ and therefore we have isomorphisms 
${}_S\unit \otimes {}_{\mathfrak q}\unit \simeq {}_S\unit $ and 
$\unit_{\mathfrak q}\otimes \unit_S \simeq \unit_{\mathfrak q}$ by 
Corollary \ref{cor:compl_comp}.
\qed

The following consequence is a local version of the more trivial remark that $K^G_*(A)\simeq 0$ for an $A\in \tcat^G$ implies $K^H_*(\Res{G}{H} A)$ $\simeq$ $ 0$. 

\begin{cor} \label{cor:central_restriction}
In the situation of Lemma \ref{lemma:central_restriction}, if $A\in \tcat^G$ and $K^G_*(A)_{\mathfrak p}\simeq 0$ then $K^H_*(\Res{G}{H}A)_{\mathfrak q}$ $\simeq$ $0$.
\end{cor}

\proof
Since $\{\unit,T(\unit)\}$ generates $\tcat^G$, $K^G_*(A)_{\mathfrak p}= K^G_*(\unit_{\mathfrak p}\otimes A)\simeq 0$ implies $\unit_{\mathfrak p}\otimes A\simeq 0$ and  therefore $ \Res{G}{H}(\unit_{\mathfrak p})\otimes \Res{G}{H}(A) \simeq 0  $. Hence, by the second isomorphism in the lemma, $\unit_{\mathfrak q}\otimes \Res{G}{H}(A)\simeq 0$ and consequently $K^H_*(\Res{G}{H}A)_{\mathfrak q}\simeq 0$.
\qed


Next, we prove $\mathfrak p$-local versions of a couple of results of \cite{meyernest-hom} which will be put to good use in the following two sections. \\

Consider the homological pair $(\tcat^G_{\mathfrak p},\mathcal I)$ with $\mathcal I:=\ker (K_*^G(?)_{\mathfrak p})$ (see Def.\ \ref{defi:hlg_alg}).
Denote by $R(G)_{\mathfrak p}\Modules^{\Z/2}_{\infty}$ the stable abelian category of $\Z/2$-graded countable (indicated by~``$\infty$'') $R(G)_{\mathfrak p}$-modules and degree-zero homomorphisms.

\begin{prop} \label{prop:univ}
The functor $h:=K^G_*(?)_{\mathfrak p}\simeq K^G_* :\tcat^G_{\mathfrak p}\to R(G)_{\mathfrak p}\Modules^{\Z/2}_{\infty}$ is the universal $\mathcal I$-exact (stable homological) functor on $\tcat^G_{\mathfrak p}$. 
Moreover, $h$ restricts to an equivalence $\mathcal P_{\mathcal I}\simeq \Projobj \big(R(G)_{\mathfrak p}\Modules_{\infty}^{\Z/2}\big)$, and, for every $A\in \tcat^G_{\mathfrak p}$, it induces a bijection between isomorphism classes of projective resolutions of $h(A)$ in $R(G)_{\mathfrak p}\Modules^{\Z/2}_{\infty}$ and isomorphism classes of $\mathcal I$-projective resolutions of $A$ in $\tcat^G_{\mathfrak p}$.
\end{prop}

\proof
We use Meyer and Nest's criterion \cite[Theorem 57]{meyernest-hom}. Since $\tcat^G_{\mathfrak p}$ is idempotent complete (having arbitrary countable coproducts); since the abelian category $R(G)_{\mathfrak p}\Modules_{\infty}^{\Z/2}$ has enough projectives (being: graded modules that are degree-wise $R(G)_{\mathfrak p}$-projective), and since $h$ is obviously an $\mathcal I$-exact stable homological functor, in order to derive the universality of $h$ from the cited theorem it remains to find for $h$ a partial left adjoint 
\begin{equation*}
h^{\dagger}\;:\; \Projobj \big(R(G)_{\mathfrak p}\Modules_{\infty}^{\Z/2}\big)\longrightarrow \tcat^G_{\mathfrak p}
\end{equation*}
defined on projective objects, such that
\begin{equation} \label{dagger}
h\circ h^{\dagger}(P)\simeq P
\end{equation}
naturally in~$P$. Since every projective in $R(G)_{\mathfrak p}\Modules_{\infty}^{\Z/2}$ is a  direct factor of a coproduct of copies of $R(G)_{\mathfrak p}(0)$ and $R(G)_{\mathfrak p}(1)$ (i.e., $R(G)_{\mathfrak p}$ concentrated in $\Z/2$-degree 0 and 1 respectively), and since $h$ preserves coproducts, it suffices to define $h^{\dagger}$ on the latter two graded modules (\cite[Remark 58]{meyernest-hom}).

Set $h^{\dagger}(R(G)_{\mathfrak p}(i)):=T^{i}(\unit_{\mathfrak p})$ for $i=0,1$, where $\unit_{\mathfrak p}\in \tcat^G_{\mathfrak p} $ is the $\mathfrak p$-localization of the tensor unit as in Theorem \ref{thm:locKth}. Then indeed, the partially defined $h^{\dagger}$ (extended to a functor in the evident way) is left adjoint to~$h$, because for all $A= \unit_{\mathfrak p}\otimes A\in \tcat^G_{\mathfrak p}$ we have
\begin{eqnarray*}
\KK^G(h^{\dagger}\big(R(G)_{\mathfrak p}(i)), A \big) 
& = &  \KK^G( T^i \unit_{\mathfrak p},\unit_{\mathfrak p}\otimes A)   \\
& \simeq & \KK^G( T^i \unit , \unit_{\mathfrak p}\otimes A) \\
& \simeq & K^G_i (A)_{\mathfrak p}  = \Hom_{R(G)} \big(R(G)(i), h(A) \big),
\end{eqnarray*}
 by Proposition \ref{prop:compl} (a) and Theorem \ref{thm:locKth}~(h).
We immediately verify \eqref{dagger}:
\begin{equation*}
hh^{\dagger}\big(R(G)_{\mathfrak p}(i)\big)= \KK^G_*(\unit, T^{i}\unit_{\mathfrak p})\simeq R(G)_{\mathfrak p}(i)
\quad \quad (i=0,1).
\end{equation*}
Thus $h$ is the universal $\mathcal I$-exact functor. 
The other claims in the proposition follow from this one, see \cite[Thm.\ 59]{meyernest-hom}.
\qed

We can use the latter proposition to compute left derived functors with respect to $\mathcal I=\ker(h)$, as follows:

\begin{prop}\label{prop:identif}
Let $F:\tcat_{\mathfrak p}^G\to \Ab$ be a homological functor which preserves small coproducts. Then for every $n\geq 0$ there is a canonical isomorphism
\begin{equation} \label{identif}
\mathsf L_n^{\mathcal I} F_*\simeq \Tor^{R(G)_{\mathfrak p}}_n\big(F_*(\unit_{\mathfrak p}), h(?)\big)
\end{equation} 
of functors $\tcat_{\mathfrak p}^G\to\Ab^{\Z/2}$. (On the left hand side we have the left derived functors of $F_*$ with respect to $\mathcal I=\ker(h)$; on the right hand side, the left derived functors of the usual tensor product  of graded modules, i.e., the homology of $\otimes^{\mathsf L}_{R(G)_{\mathfrak p}}$; the $R(G)_{\mathfrak p}$-action on $F_*(\unit_{\mathfrak p})$ is  induced by the functoriality of~$F$, cf.\ Rem.\ \ref{remark:iso-of-modules}.) 
\end{prop}

\proof (Note by inspecting the definitions that $\mathsf L^{\mathcal I}_n(F_*)=(\mathsf L^{\mathcal I}_nF)_*$.)
We have proved above that $h$ is the universal $\mathcal I$-exact functor. It follows that every homological functor $F:\tcat_{\mathfrak p}^G\to \mathcal A$ extends (up to isomorphism, uniquely) to a right exact functor
$$\tilde F\;:\; R(G)_{\mathfrak p}\Modules_{\infty}^{\Z/2}\longrightarrow \mathcal A$$
such that $\tilde F\circ h(P)= F(P)$ for all $\mathcal I$-projective objects~$P$; this functor $\tilde F$ is stable, resp.\ commutes with coproducts, if so does~$F$.
Moreover, there are canonical isomorphisms
\begin{equation} \label{derfunctorsF}
\mathsf L_n^{\mathcal I} F_*\simeq  (\mathsf L_n\tilde F_*)\circ h
\end{equation}
for all $n\in\Z$. (See \cite[Theorem 59]{meyernest-hom} for these results).
Therefore we are left with computing $\tilde F_*$ and its left derived functors, in the case where $\mathcal A$ is the category of abelian groups. 

\begin{lemma}\label{lemma:modules} There is a natural isomorphim
\begin{equation}\label{formulaHbar}
\tilde F_*(M)\simeq F_*(\unit_{\mathfrak p})\otimes_{R(G)_{\mathfrak p}}M
\end{equation}
of graded abelian groups, for $M\in R(G)_{\mathfrak p}\Modules^{\Z/2}_{\infty}$.
\end{lemma}
To prove the lemma, notice first that \eqref{formulaHbar} holds for the free module $M=R(G)_{\mathfrak p}$ (set in degree zero), because there are canonical isomorphisms of graded $R(G)_{\mathfrak p}$-modules
$$\tilde F_*(R(G)_{\mathfrak p})
= \tilde F_* \circ h (\unit_{\mathfrak p}) 
= F_*(\unit_{\mathfrak p})
\simeq F_*(\unit_{\mathfrak p})\otimes_{R(G)_{\mathfrak p}} R(G)_{\mathfrak p} .$$
We may extend this to all $\Z/2$-graded free modules in the evident way. Since both $\tilde F_*$ and $F_*(\unit_{\mathfrak p})\otimes(?)$ are right exact functors, we can compute them -- and we can extend the natural isomorphism \eqref{formulaHbar} -- for general graded modules $M$ by using free presentations $P\to P'\to M\to 0$. \qed

Proposition \ref{prop:identif} follows now from Lemma \ref{lemma:modules}:  
by taking left derived functors of \eqref{formulaHbar} we get $\mathsf L_n \tilde F_*\simeq \Tor_n^{R(G)_{\mathfrak p}}(F_*(\unit_{\mathfrak p}),\,?)$, and by combining this with \eqref{derfunctorsF} we find the predicted isomorphism \eqref{identif}.
\qed

\begin{remark} \label{remark:iso-of-modules}
Let $F:\tcat_{\mathfrak p}^G\to \Ab$ be an additive functor. Since $\tcat_{\mathfrak p}^G$ is an $R(G)_{\mathfrak p}$-linear category, $F$ lifts to $R(G)_{\mathfrak p}\Modules^{\Z/2}$, simply via $r\cdot a:=F(r\cdot \id_A)(a)$ for all $r\in R(G)_{\mathfrak p}$ and $a\in F(A)$. This is for instance how we regard $F_*(\unit_{\mathfrak p})$ as a graded $R(G)_{\mathfrak p}$-module in Proposition \ref{prop:identif}.
It is clear from the proof that the isomorphism \eqref{identif} is actually an isomorphism of graded $R(G)_{\mathfrak p}$-modules.
\end{remark}

The same arguments provide an analog statement for contravariant functors. We leave the details of the proof to the reader (cf.\ \cite[Thm.\ 72]{meyernest-hom}):

\begin{prop} \label{prop:identif_contrav}
Let $F: (\tcat_{\mathfrak p}^G)^{\op}\to \Ab$ be a homological functor sending small coproducts in $\tcat^G_{\mathfrak p}$ to products. Then for every $n\geq 0$ there is an isomorphism
\begin{equation*}
\mathsf R^n_{\mathcal I} F_* \simeq \Ext^n_{R(G)_{\mathfrak p}} \big( h( \textrm{\emph{?`}} ), F_*(\unit)\big)
\end{equation*} 
of contravariant functors from $\tcat_{\mathfrak p}^G$ to $\Z/2$-graded $R(G)_{\mathfrak p}$-modules. (The graded Ext on the right are the derived functors of the graded Hom $\Hom^*_{R(G)_{\mathfrak p}}(\textrm{\emph{?`}}, F_*(\unit))$.)
\qed
\end{prop}

\subsection{The Phillips-K\"unneth formula}
\label{subec:PKF}

We derive from the above theory a new version of a theorem of N.\,C.\ Phillips (\cite[Theorem 6.4.6]{phillips-free}). Our theorem and that of Phillips differ only in the technical assumptions on the $\Cstar$-algebras involved;  we don't know how these compare precisely, but we suspect that neither set of hypotheses implies the other. 

Phillips' theorem is about the following data, whose relevance will be explained at the beginning of~$\S$\ref{subsec:supp_datum}.

\begin{defi} \label{defi:local_pair}
A \emph{local pair} $(S,\mathfrak q)$ consists of a finite cyclic group $S$ and  a prime ideal $\mathfrak q\in \spec(R(S))$ such that, if $S'\leq S$ is a subgroup with the property that $(\Res{S}{S'})^{-1}(\mathfrak q')=\mathfrak q$ for some $\mathfrak q'\in \spec(R(S'))$, then $S'=S$. (Here $\Res{S}{S'}:R(S)\to R(S')$ is the usual restriction ring homomorphism; of course, it coincides with the functor $\Res{S}{S'}: \KK^S\to \KK^{S'}$ at $R(S)=\KK^S(\unit,\unit)$.)  
\end{defi}

\begin{lemma} \label{lemma:loc_ring_hered}
Let $(S,\mathfrak q)$ be a local pair. Then the local ring $R(S)_{\mathfrak q}$ is a discrete valuation ring or a field; in particular, it is hereditary (that is,  every submodule of a projective $R(S)_{\mathfrak q}$-module is again projective).
\end{lemma}
\proof
See \cite[Prop.\ 6.2.2]{phillips-free}, where it is proved that, under the above hypothesis, $R(S)_{\mathfrak q}$ is isomorphic to the localization at a prime ideal of $\Z[\zeta]$, the subring  of $\C$ generated by a primitive $n$th root of unity~$\zeta$, where $n=|S|$. The claims follow because $\Z[\zeta]$ is a Dedekind domain (cf.\ \cite[Lemma 6.4.2]{phillips-free}).
\qed

\begin{thm} \label{thm:phku}
\emph{(Phillips-K\"unneth Formula).}
Let $(S,\mathfrak q)$ be a local pair. Then for all $A\in \mathcal T^S$ and $B\in \KK^S$ there is a natural short exact sequence
\begin{equation*}
\xymatrix{
K^S_*(A)_{\mathfrak q}\otimes_{R(S)_{\mathfrak q}}K^S_*(B)_{\mathfrak q} \ar@{ >->}[r] &
  K^S_*(A\otimes B)_{\mathfrak q} \ar@{->>}[r]^-{+1} &
    \Tor_1^{R(S)_{\mathfrak q}}(K^S_*(A)_{\mathfrak q}, K^S_*(B)_{\mathfrak q})
}
\end{equation*}
of $\Z/2$-graded $R(S)_{\mathfrak q}$-modules which splits unnaturally (the $+1$ indicates a map of $\Z/2$-degree one).
\end{thm}

\begin{lemma} \label{lemma:special}
It suffices to prove the theorem for the special case $A,B \in \tcat^S_{\mathfrak q}$. 
\end{lemma}

\proof
Let $A\in \tcat^S$ and $B\in \KK^S$. 
Let $LB \to B \to RB\to TLB$ be the natural distinguished triangle with $LB\in \tcat^S$ and $K^S_*(RB)\simeq 0$ (Thm.\ \ref{thm:compl_Ktheory}).
Since $LB\to B$ induces an isomorphism $K^S_*(LB)\simeq K^S_*(B)$, we may substitute $LB$ for $B$ in the first and third terms of the sequence.  Note that the subcategory $\{X\in \KK^S \mid K^S_*(X\otimes RB)\simeq 0\}$ is localizing and contains $\unit$, hence it contains $\tcat^S$.
Therefore $LB\to B$ also induces an isomorphism $K^S_*(A\otimes LB)\simeq K^S_*(A\otimes B)$.
Hence it suffices to prove the existence and split exactness of the sequence for $A,B\in \tcat^S$.

Now, if $A,B\in \tcat^S$ then $K^S_*(\unit_{\mathfrak q}\otimes A)_{\mathfrak q}= K^S_*(A)_{\mathfrak q}$, $K^S_*(\unit_{\mathfrak q}\otimes B)_{\mathfrak q}= K^S_*(B)_{\mathfrak q}$ and $K^S_*(\unit_{\mathfrak q}\otimes A\otimes \unit_{\mathfrak q} \otimes B)_{\mathfrak q}= K^S_*(A\otimes B)_{\mathfrak q}$ by Theorem \ref{thm:locKth}, so we may as well substitute $\unit_{\mathfrak q}\otimes A\in \tcat^S_{\mathfrak q}$ for $A$ and $\unit_{\mathfrak q}\otimes B\in \tcat^S_{\mathfrak q}$ for $B$.
\qed


\proof[Proof of Theorem \ref{thm:phku}]
By the previous lemma we can assume that $A\in \tcat^S_{\mathfrak q}$.
We wish to apply Theorem \ref{thm:hereditary} (a) to the homological pair $(\tcat^S_{\mathfrak q},\mathcal I:= \ker(K^S_*(?)_{\mathfrak q}))$ and the homological functor $F:= K^S_*(?\otimes B)_{\mathfrak q}$.

By Prop.\ \ref{prop:univ},  $h:=K^S_*(?)_{\mathfrak q}:\tcat^S_{\mathfrak q}\to R(S)_{\mathfrak q}\Modules^{\Z/2}_{\infty}$ is the universal $\mathcal I$-exact functor and therefore it induces a bijection between isomorphism classes of projective resolutions of the graded $R(S)_{\mathfrak q}$-module $K^S_*(A)_{\mathfrak q}$ and isomorphism classes of $\mathcal I$-projective resolutions of $A$.
By Lemma \ref{lemma:loc_ring_hered} every $R(S)_{\mathfrak q}$-module has a projective resolution of length one, so $A$ has an $\mathcal I$-projective resolution of length one. Since $A\in \tcat^S_{\mathfrak p}=\langle \unit_{\mathfrak q}\rangle_{\loc}=\langle \mathcal P_{\mathcal I}\rangle_{\loc}$, it satisfies the hypothesis of Theorem \ref{thm:hereditary}.
Therefore there exists a natural short exact sequence $0\to \mathsf L^{\mathcal I}_0F(A)\to F(A) \to \mathsf L^{\mathcal I}_1F(TA)\to 0$. It remains to identify the derived functors of $F=K^S_*(?\otimes B)_{\mathfrak q}$ and to show that the sequence splits. 
According to Proposition \ref{prop:identif} (applied to the homological functor $K^S_0(?\otimes B)_{\mathfrak q}$), we have a natural isomorphism
\begin{eqnarray*}
\mathsf L_i^{\mathcal I}F(A) 
&\simeq & \Tor^{R(S)_{\mathfrak q}}_i(K^S_*( \unit_{\mathfrak q}\otimes B)_{\mathfrak q}, h_*(A)) \\
&= & \Tor^{R(S)_{\mathfrak q}}_i (K^S_{*}(B)_{\mathfrak q}, K^S_*(A)_{\mathfrak q}) \\
&=& \Tor^{R(S)_{\mathfrak q}}_i(K^S_*(A)_{\mathfrak q},K^S_{*}(B)_{\mathfrak q})
\end{eqnarray*}
of graded $R(S)_{\mathfrak q}$-modules for $i=0,1$,  as claimed. As for the splitting, we can use the same argument as in \cite[$\S$23.11]{blackadar}. We postpone this to Corollary \ref{cor:splitting}, which requires the (unsplit) universal coefficient theorem. \qed

\begin{thm}[Universal Coefficient Theorem, UCT] \label{thm:UCT}
Let $(S,\mathfrak q)$ be a local pair. For every $A\in \tcat^S$ and $B\in \KK^S$ 
there exists a natural short exact sequence
\begin{equation*}
\xymatrix@1{
 \Ext^1_{R(S)_{\mathfrak q}}(K^S_*(A)_{\mathfrak q}, K^S_*(B)_{\mathfrak q}) \ar@{ >->}[r]^-{+1} &
  \KK^S_*(A,B)_{\mathfrak q} \ar@{->>}[r] &
   \Hom^*_{R(S)_{\mathfrak q}}(K^S_*(A)_{\mathfrak q},K^S_*(B)_{\mathfrak q})
}
\end{equation*}
of $\Z/2$-graded $R(S)_{\mathfrak q}$-modules.
\end{thm}

\proof
The proof  is quite similar to that of Theorem \ref{thm:phku}. Just as before in Lemma \ref{lemma:special} we reduce to the case $A,B \in \tcat^S_{\mathfrak q}$, but then we use Theorem \ref{thm:hereditary} (c) (for both $B$ and $TB$) to produce the short exact sequence and Proposition \ref{prop:identif_contrav} to identify its right and left terms as required (cf.\ \cite[Thm.\ 72]{meyernest-hom}).
\qed

The UCT has corollaries familiar from ordinary $K$-theory (cf.\ \cite[$\S$23]{blackadar}). 

\begin{cor}\label{cor:realization}
Let $M$ be any countable $\Z/2$-graded $R(S)_{\mathfrak q}$-module. Then there exists an object $A\in \tcat^S_{\mathfrak q}$ such that $K^S_*(A)=K^S_*(A)_{\mathfrak q}\simeq M$.
\end{cor}

\proof
Consider a projective (i.e., free) resolution $0\to Q\to P \to M\to 0$ in $R(S)_{\mathfrak q}\Modules^{\Z/2}_{\infty}$. 
Applying $h^{\dagger}$ (see the proof of Proposition \ref{prop:univ}) we obtain a morphism $f:h^{\dagger}Q\to h^{\dagger}P$ between $\mathcal I$-projective objects in $\tcat^S_{\mathfrak q}$. Now apply $h=K^S_*(?)_{\mathfrak q}$ to the distinguished triangle $h^{\dagger}Q\to h^{\dagger}P\to \cone(f) \to Th^{\dagger}Q $ to get the exact sequence
$
Q \to P \to K^S_*(\cone(f))_{\mathfrak q} \to Q[1] \to P[1].
$ The rightmost map is injective and therefore $K^S_*(\cone(f))_{\mathfrak q}\simeq M$.
\qed

\begin{cor}\label{cor:iso}
Consider objects $A,B\in \tcat^S_{\mathfrak q}$ such that $K^S_*(A)_{\mathfrak q}\simeq K^S_*(B)_{\mathfrak q}$. Then there exists an isomorphism $A\simeq B$ in $\tcat^S_{\mathfrak q}$.
\end{cor}

\proof
Because of the surjectivity of the second homomorphism in the UCT (in degree zero), we may lift the isomorphism $K^S_*(A)_{\mathfrak q}\simeq K^S_*(B)_{\mathfrak q}$ to a map $f:A\to B$ in~$\tcat^S_{\mathfrak q}$. Since $\{\unit, T(\unit)\}$ generates $\tcat^S$, the condition 
$\cone(f)\simeq 0$ is equivalent to 
$\KK^S_*(\unit, \cone(f))=K^S_*(\cone(f))_{\mathfrak q}\simeq 0$. But $K^S_*(f)_{\mathfrak q}$ is an isomorphism by construction, hence $f:A\simeq B$.
\qed

\begin{cor} \label{cor:decomp} 
Let $A\in \tcat^S_{\mathfrak q}$, and assume that there is an isomorphism $K_*^S(A)_{\mathfrak q}\simeq M_1\oplus M_2$ of graded $R(S)_{\mathfrak q}$-modules. Then there exists in $\tcat^S_{\mathfrak q}$ a decomposition $A\simeq A_1\oplus A_2$ with $K_*^S(A_i)_{\mathfrak q}\simeq M_i$ ($i=1,2$).
\end{cor}

\proof
Use Corollary \ref{cor:realization} to get $A_i\in \tcat^S_{\mathfrak q}$ with $K_*^S(A_i)\simeq M_i$ ($i=1,2$). Now employ Corollary \ref{cor:iso}.
\qed

\begin{cor} \label{cor:splitting}
The short exact sequences in the Phillips-K\"unneth Theorem \ref{thm:phku} and the Universal Coefficient Theorem \ref{thm:UCT} are (unnaturally) split.
\end{cor}

\proof
If $\tilde A\in \tcat^S_{\mathfrak q}$, according to Corollary \ref{cor:decomp} the degree-wise decomposition $K^S_*(\tilde A)_{\mathfrak q}= K^S_0(\tilde A)_{\mathfrak q}(0) \oplus K^S_1(\tilde A)_{\mathfrak q}(1) $
 can be realized by a decomposition 
 $\tilde A\simeq A_0\oplus A_1$ in $\tcat^S_{\mathfrak q}$. Let $A\in \tcat^S$. Now we apply the preceding to 
 $\tilde A:=\unit_{\mathfrak q}\otimes A\in \tcat^S_{\mathfrak q}$
  and appeal to Remark \ref{rem:splitting}.
\qed

\subsection{The residue field object at a prime ideal}
\label{subsec:res_field_object}

Fix a local pair $(S,\mathfrak q)$, as in Def.\ \ref{defi:local_pair}. That is: $S$ is a cyclic group and $\mathfrak q\in\spec R(S)$ does not lie above any $\mathfrak q'\in \spec R(S')$ with $S'< S$ a proper subgroup. 
Denote by $k(\mathfrak q):=R(S)_{\mathfrak q}/\mathfrak qR(S)_{\mathfrak q}$ the residue field of $R(S)$ at the prime ideal~$\mathfrak q$. The following lemma is an immediate consequence of Corollary \ref{cor:realization}. Together with the Phillips-K\"unneth formula, it is the key ingredient needed for the construction of the support $\sigma_G$ in Theorem \ref{thm:gen_supp_intro}.

\begin{lemma} \label{lemma:res_field_obj}
There exists an object $\kappa_{\mathfrak q}\in \tcat^S_{\mathfrak q}$ with the property that $K^S_0(\kappa_{\mathfrak q})\simeq k(\mathfrak q)$ and $K^S_1(\kappa_{\mathfrak q})\simeq0$.
\qed
\end{lemma}

\begin{defi}\label{defi:rfo}
We call such an object $\kappa_{\mathfrak q}$ a \emph{residue field object at $(S,\mathfrak q)$}. By Corollary \ref{cor:iso}, it is uniquely determined by $(S,\mathfrak q)$ up to isomorphism.
\end{defi}

\begin{prop} \label{prop:field_object}
For every $A\in \tcat^S$, the product $\kappa_{\mathfrak q}\otimes A$ is isomorphic in $\tcat^S$ to a countable coproduct of translated copies of $\kappa_{\mathfrak q}$.
\end{prop}

\proof
Note that $\kappa_{\mathfrak q}\otimes A \in \tcat^S_{\mathfrak q}$.
Applied to the objects $\kappa_{\mathfrak q}$ and $A$, the Phillips-K\"unneth split short exact sequence (Thm.\ \ref{thm:phku}) implies that the $\Z/2$-graded $R(S)_{\mathfrak q}$-module $K^S_*(\kappa_{\mathfrak q}\otimes A)$ is isomorphic to a $\Z/2$-graded $k(\mathfrak q)$-vector space, which has the form $\coprod_{I_0} k(\mathfrak q)(0)\oplus \coprod_{I_1} k(\mathfrak q)(1)$ for some countable index sets $I_0$ and $I_1$. The latter vector space can be realized in $\tcat^S_{\mathfrak q}$ as the object $B:=\coprod_{I_0} \kappa_{\mathfrak q}\oplus \coprod_{I_1} T(\kappa_{\mathfrak q})$. Since $\kappa_{\mathfrak q}\otimes A$ and $B$ both lie in $\tcat^S_{\mathfrak q}$ and have isomorphic $K$-theory, by Corollary \ref{cor:iso} of the UCT they must be isomorphic. 
\qed

\begin{prop} \label{prop:prod_formula}
Let $(S,\mathfrak q)$ be a local pair. Then for every two objects $A,B\in \mathcal T^S$ there exists a (non natural) isomorphism
\begin{equation*}
K^S_*(\kappa_{\mathfrak q}\otimes A \otimes B) \;\simeq \; K^S_*(\kappa_{\mathfrak q}\otimes A) \;\hat \otimes \; K^S_*(\kappa_{\mathfrak q}\otimes B)
\end{equation*}
of $\Z/2$-graded $k(\mathfrak q)$-vector spaces. Here $\hat \otimes$ denotes the usual tensor product of graded vector spaces, given by $(V\hat \otimes W)_{\ell}= \bigoplus_{i+j=\ell}V_i\otimes_{k(\mathfrak q)} V_j$.
\end{prop}

\proof
To simplify notation, we write $\kappa:=\kappa_{\mathfrak q}$ and $k:=k(\mathfrak q)$. Choose isomorphisms
\begin{equation*}
\kappa\otimes A \simeq \coprod_{n_0} \kappa \oplus \coprod_{n_1} T(\kappa) \quad \textrm{ and }\quad
\kappa\otimes B \simeq \coprod_{m_0} \kappa \oplus \coprod_{m_1} T(\kappa)
\end{equation*}
in $\tcat^S$ as provided by Proposition \ref{prop:field_object}. Then
\begin{eqnarray*}
\kappa\otimes A \otimes B
 &\simeq&  \big( \coprod_{n_0} \kappa \oplus \coprod_{n_1} T(\kappa) \big) \otimes B \\
&\simeq& \big( \coprod_{n_0} \kappa\otimes B\big) \oplus \big( \coprod_{n_1} T(\kappa\otimes B) \big) \\
&\simeq & \coprod_{n_0} \big(  \coprod_{m_0} \kappa \oplus \coprod_{m_1} T(\kappa)\big) \oplus \coprod_{n_1} \big(  \coprod_{m_0} T(\kappa) \oplus \coprod_{m_1} \kappa\big) \\
&\simeq & \coprod_{n_0m_0 + n_1m_1} \kappa \quad \oplus \coprod_{n_0m_1 + n_1m_0} T(\kappa).
\end{eqnarray*}
Since $K^S_*(\kappa)\simeq k(0)$ and $K^S_*(T\kappa)\simeq k(1)$ (where, as before, $V(i)$ stands for the $k$-vector space $V$ set in degree $i\in \Z/2$), we obtain
\begin{equation*}
K^S_*(\kappa\otimes A\otimes B) \simeq 
\coprod_{n_0m_0 + n_1m_1} k(0) \quad \oplus \coprod_{n_0m_1 + n_1m_0} k(1).
\end{equation*}
The right hand side of the equation is computed similarly: 
\begin{eqnarray*}
K^S_*(\kappa\otimes A)\; \hat \otimes \; K^S_*(\kappa\otimes B)  &\simeq& 
\big(\coprod_{n_0} k(0) \oplus \coprod_{n_1} k(1)\big)\; \hat \otimes \;  
\big(  \coprod_{m_0} k(0) \oplus \coprod_{m_1} k(1) \big)
\\
&\simeq &
\coprod_{n_0m_0 + n_1m_1} k(0) \quad \oplus \coprod_{n_0m_1 + n_1m_0} k(1)
\end{eqnarray*}
using that $k(i)\;\hat \otimes\; k(j)\simeq k(i+j)$. We see that the two sides are isomorphic.
\qed

We also record the following consequence of the Phillips-K\"unneth theorem.

\begin{cor} \label{cor:kkey}
 Let $A\in \tcat^S$. Then $K^S_*(\kappa_{\mathfrak q}\otimes A)\simeq 0$  if and only if  the derived tensor product $k(\mathfrak q)\otimes^{\mathsf L}_{R(S)_{\mathfrak q}}K^S_*(A)_{\mathfrak q}= k(\mathfrak q)\otimes^{\mathsf L}_{R(S)}K^S_*(A)$ is zero.

\end{cor}
\proof
Since $\kappa_{\mathfrak q}\simeq \unit_{\mathfrak q}\otimes \kappa_{\mathfrak q}$, we may substitute $A$ with $\unit_{\mathfrak q}\otimes A$ and $K^S_*(\kappa_{\mathfrak q}\otimes A)$ with $K^S_*(\kappa_{\mathfrak q}\otimes A)_{\mathfrak q}$. 
By the Phillips-K\"unneth formula \ref{thm:phku}, $K^S_*(\kappa_{\mathfrak q}\otimes A)_{\mathfrak q}$ vanishes if and only if $\Tor_i^{R(S)_{\mathfrak q}}(k(\mathfrak q),K^S_*(A)_{\mathfrak q})\simeq 0 \quad (i=0,1)$. The latter Tor modules are by definition the homology of the complex $k(\mathfrak q)\otimes^{\mathsf L}_{R(S)_{\mathfrak q}}K^S_*(A)_{\mathfrak q}$. 
\qed\newline

\section{First results for finite groups}
\label{sec:first_results}

\subsection{The nice support $(\spec R(G),\sigma_G)$ on~$\tcat^G$}
\label{subsec:supp_datum}

We are now ready to prove Theorem \ref{thm:gen_supp_intro} of the introduction. 
We fix an arbitrary \emph{finite} group~$G$ and consider the compactly generated $\otimes$-triangulated category $\tcat^G=\langle \unit \rangle_{\loc}\subset \KK^G$ of $\S$\ref{subsec:cats}. 

In \cite{segal-spec}, it is shown that for every prime ideal $\mathfrak p\in \spec (R(G))$ there exists a cyclic subgroup $S\leq G$, unique up to conjugacy in~$G$ (let us call it the \emph{source}\footnote{In \emph{loc.\ cit.} Segal calls it the \emph{support} of $\mathfrak p$, but surely the reader of this article will forgive us for avoiding charging this poor word with yet another meaning.} of~$\mathfrak p$), such that: There exists a prime ideal $\mathfrak q\in\spec (R(S))$ with  $(\Res{G}{S})^{-1}(\mathfrak q)=\mathfrak p$, and moreover $S$ is minimal (with respect to inclusion) among the subgroups of $G$ with this property. It follows that $\mathfrak q$ also cannot come from any proper subgroups of~$S$, i.e., the source of such a $\mathfrak q\in \spec(R(S))$ is $S$ itself. 

\begin{notation} In the following, for a $\mathfrak p\in \spec(R(G))$ and a fixed cyclic subgroup $S=S(\mathfrak p)$ of $G$ in the conjugacy class of the source of~$\mathfrak p$, we shall denote by
\begin{equation*}
\local(\mathfrak p):= \{\mathfrak q\in \spec(R(S(\mathfrak p)))\mid (\Res{G}{S(\mathfrak p)})^{-1}(\mathfrak q)=\mathfrak p\}
\end{equation*}
 the fiber in $\spec(R(S(\mathfrak p)))$ over the point $\mathfrak p\in \spec(R(G))$.
\end{notation}

Note that the pair $(S(\mathfrak p),\mathfrak q)$, for  any $\mathfrak q\in \local(\mathfrak p)$, is a local pair as in Definition \ref{defi:local_pair}. In particular, we can apply to it all the results of $\S$\ref{subsec:res_field_object}, such as the existence of a residue field object $\kappa_{\mathfrak q}\in \tcat^{S(\mathfrak p)}_{\mathfrak q}$ (Lemma \ref{lemma:res_field_obj}).

\begin{defi} \label{defi:sigma}
For a local pair $(S,\mathfrak q)$, denote by $\mathcal A(S,\mathfrak q)$ the stable abelian category of countable $\Z/2$-graded $k(\mathfrak q)$-vector spaces. Write
\begin{equation*}
F_{(S,\mathfrak q)} : \tcat^S \longrightarrow \mathcal A(S,\mathfrak q)
\end{equation*}
for the stable homological functor sending $B\in \tcat^S$ to $K^S_*(\kappa_{\mathfrak q}\otimes B)$. Now for every $\mathfrak p\in \spec(R(G))$, choose a $\mathfrak q= \mathfrak q(\mathfrak p) \in \local(\mathfrak p)$ and consider the functor 
\begin{equation*} \label{F_p}
 F_{\mathfrak p} := F_{(S(\mathfrak p),\mathfrak q(\mathfrak p))}\circ \Res{G}{S(\mathfrak p)}
\quad :\quad
\xymatrix@1{
\tcat^G \ar[r] &  \mathcal A (S(\mathfrak p), \mathfrak q(\mathfrak p)) =:\mathcal A(\mathfrak p).
}
\end{equation*}
Finally, define the support $\sigma_G$ by 
\begin{eqnarray*}
\sigma_G(A)&:=& \{\mathfrak p\mid F_{\mathfrak p}(A)\not\simeq 0\} \\
&=& \{\mathfrak p\mid K_*^{S(\mathfrak p)}\big(\kappa_{ \mathfrak q(\mathfrak p)}\otimes \Res{G}{S(\mathfrak p)}A\big)\not\simeq 0\} \\
&=& \{\mathfrak p\mid \kappa_{\mathfrak q(\mathfrak p)}\otimes \Res{G}{S(\mathfrak p)}(A)\not\simeq 0\}
\;\;\; \subset\;\;\; \spec (R(G))
\end{eqnarray*}
for every object $A\in\tcat^G$. 
\end{defi}

\begin{remark} \label{rem:indep}
The set $\sigma_G(A)\subset \spec(R(G))$ only depends on the group $G$ and the object $A\in \tcat^G$, not on the choices of $S(\mathfrak p)$, $\mathfrak q(\mathfrak p)\in \local(\mathfrak p)$ or $\kappa_{\mathfrak q(\mathfrak p)}$. By Cor.\ \ref{cor:kkey}, for fixed $(S,\mathfrak q)=(S(\mathfrak p),\mathfrak q(\mathfrak p))$ the vanishing of $F_{\mathfrak p}(A)$ only depends on the $R(S)$-module $K^S_*(\kappa_{\mathfrak q})\simeq k(\mathfrak q)$, not on the choice of $\kappa_{\mathfrak q}\in \tcat^G_{\mathfrak q}$.
Now let $(S,\mathfrak q)$ and $(S',\mathfrak q')$ be two choices. As we already noted, if $S$ and $S'$ are two cyclic subgroups of $G$, both representing the source of $\mathfrak p$, then $S$ and $S'$ are conjugate in~$G$; moreover, any two primes $\mathfrak q_1,\mathfrak q_2\subset \spec(R(S))$ lying above $\mathfrak p$  are also conjugate by the induced action of some element of the normalizer $N_G(S)$  (\cite[Prop.\ 3.5]{segal-spec}). Combining the two, we easily find an isomorphism  
 $\phi:S \stackrel{\sim}{\to}S'$, $s\mapsto g^{-1}sg$ inducing a $\otimes$-triangulated isomorphism $\phi^*:\KK^{S'}\simeq \KK^{S}$ such that $\phi^*\circ\Res{G}{S'} \simeq \Res{G}{S}$ and $\phi^*(\kappa_{\mathfrak q'})\simeq \kappa_{\mathfrak q}$. This shows that $\sigma_G(A)$ is independent of all choices. 
\end{remark}

\begin{thm} \label{thm:sigma}
The pair $(\spec R(G),\sigma_G)$ defines a support on $\tcat^G$ enjoying all the properties stated in Theorem \ref{thm:gen_supp_intro}. These are (S0)-(S7) of Theorem \ref{thm:abstract}, where moreover (S5) holds for any two objects:
\begin{equation*}
\sigma_G(A\otimes B) = \sigma_G(A) \cap \sigma_G(B) 
\end{equation*}
for all $A,B\in \tcat^G$. In particular, the restriction $(\spec(R(G)), \sigma_G|_{\mathcal K^G})$ defines a support datum on the subcategory $\mathcal  K^G= (\tcat^G)_c$ of compact objects.
\end{thm}

\proof
By definition, $\sigma_G$ is the support $\sigma_{\mathcal F(G)}$ induced, as in Lemma \ref{lemma:hml_supp}, by the family of functors $\mathcal F(G):=\{F_{\mathfrak p}\}_{\mathfrak p\in \spec R(G)}$. Every $F_{\mathfrak p}:\tcat^G\to \mathcal A(\mathfrak p)$ is a stable homological functor commuting with coproducts, because it is by definition a composition of a triangulated functor followed by a stable homological one, both of which preserve small coproducts. 
Thus, by Lemma \ref{lemma:hml_supp}, $\sigma_G$ satisfies properties (S0), (S2)-(S4) and (S6). Since $F_{\mathfrak p}(\unit)=k(\mathfrak q(\mathfrak p))\not\simeq 0$, (S1) holds as well.
 Moreover, every $\mathcal A(\mathfrak p)$ can be equipped with the tensor product $\hat \otimes$ of graded vector spaces, and clearly a product $V\,\hat \otimes\, W$ in $\mathcal A(\mathfrak p)$ is zero if and only if one of the two factors already is (consider bases). For any two objects $A,B\in \tcat^G$, there exists an isomorphism
\begin{eqnarray*}
F_{\mathfrak p}(A\otimes B)\simeq F_{\mathfrak p}(A)\,\hat\otimes\, F_{\mathfrak p}(B)  
\end{eqnarray*}
because of Proposition \ref{prop:prod_formula} and because restriction $\Res{G}{S(\mathfrak p)}$ is a $\otimes$-functor. It follows that $\sigma_G$ enjoys (S5) for any two objects.

It remains only to verify property (S7). We will do so in a series of lemmas.

\begin{lemma} \label{lemma:fg}
If $H$ is a finite (or compact Lie) group and $A\in\tcat^H_c$, then the $R(H)$-module $K^H_*(A)$ is finitely generated.
\end{lemma}
\proof
The proof is a routine induction on the length of the object $A\in \tcat^H_c=\langle \unit \rangle$, using that $R(H)$ is noetherian. We leave it to the reader. 
\qed

\begin{lemma} \label{lemma:sigma-cpt}
For every compact object $A\in\tcat^G_c$, we have
\begin{eqnarray*}
\sigma_G(A)
&=& \{\mathfrak p \in \spec(R(G)) \mid  K^{S(\mathfrak p)}_*(\Res{G}{S(\mathfrak p)} A)_{\mathfrak q(\mathfrak p)}\not\simeq 0 \}.
\end{eqnarray*}
\end{lemma}

\proof 
Write $S=S(\mathfrak p)$ and $\mathfrak q= \mathfrak q(\mathfrak p)$.
We know by Corollary \ref{cor:kkey} that
$F_{\mathfrak p}(A)=K^S_*(\kappa_{\mathfrak q}\otimes \Res{}{}A)\simeq 0$ is equivalent to the vanishing of $X_{\bullet}:=k(\mathfrak q)\otimes_{R(S)_{\mathfrak q}}^{\mathsf L}K^S_*(\Res{}{}A)_{\mathfrak q}$. Let us show that the latter is equivalent to $K^S_*(\Res{}{} A)_{\mathfrak q}\simeq 0$.  
Since $A$ is compact in $\tcat^G$, $\Res{}{} A$ is compact in $\tcat^S$ and therefore the $R(S)_{\mathfrak q}$-module $M:=K^S_*(\Res{}{}A)_{\mathfrak q}$ is finitely generated, by Lemma \ref{lemma:fg}. 
Since $R(S)_{\mathfrak q}$ is a noetherian ring of global dimension one (Lemma \ref{lemma:loc_ring_hered}), we find a length-one resolution of $M$ by finitely generated projectives, say $P_{\bullet}=(\cdots0\to P_1 \stackrel{d}{\to} P_0 \to 0\cdots)$. 
Moreover, since $R(S)_{\mathfrak q}$ is local and the $P_i$ finitely generated, we may choose the complex $P_{\bullet}$ to be \emph{minimal}, that is, such that $d(P_1)\subset \mathfrak mP_0$ where $\mathfrak m:=\mathfrak qR(S)_{\mathfrak q}$ denotes the maximal ideal (see \cite{roberts}).
 Now $X_{\bullet}=k(\mathfrak q)\otimes^{\mathsf L}M = k(\mathfrak q)\otimes P_{\bullet} = ( P_1/\mathfrak mP_1 \stackrel{0}{\to} P_0/\mathfrak mP_0)$; so $X_{\bullet}\simeq 0$ iff $P_i/\mathfrak mP_i=0$ $(i=0,1)$. By Nakayama (or simply because the modules $P_i$ are free), the latter condition is equivalent to $P_i\simeq 0$ $(i=0,1)$, i.e., to $M\simeq 0$.
 \qed

Finally, let us prove the remaining claim of Theorem \ref{thm:sigma}.

\begin{lemma} \label{lemma:genuine-supp}
The support $(\spec(R(G)),\sigma_G)$ satisfies (S7): for every $A\in \tcat^G_c$, the set $\sigma_G(A)$ is closed in $\spec(R(G))$.
\end{lemma}

\proof
Let $A$ be a compact object of $\tcat^G$. By Lemma \ref{lemma:sigma-cpt}, we can express the complement of $\sigma_G(A)$ as follows:
\begin{equation*}
\spec(R(G))\smallsetminus \sigma_G(A)  
= \{\mathfrak p\in \spec(R(G))\mid 
K^{S(\mathfrak p)}_*(\Res{G}{S(\mathfrak p)} A)_{\mathfrak q(\mathfrak p)}\simeq 0
\}.
\end{equation*}
Note that,  whenever $S$ is a cyclic subgroup of $G$ containing $S(\mathfrak p)$ and $\mathfrak r$ is a prime ideal in $R(S)$ such that $\mathfrak r=\Res{}{}^{-1}(\mathfrak q)$ and $\mathfrak p=\Res{}{}^{-1}(\mathfrak r)$, then 
\begin{eqnarray*}
K^{S}_*(\Res{G}{S}  A )_{\mathfrak r }\simeq 0
 \quad \Longrightarrow \quad 
K^{S(\mathfrak p)}_*(\Res{G}{S(\mathfrak p)}  A )_{\mathfrak q } \simeq 0
\end{eqnarray*}
by Corollary \ref{cor:central_restriction}.
Hence, by the minimality and uniqueness, up to conjugacy in $G$, of the pair $(S(\mathfrak p), \mathfrak q(\mathfrak p))$ (see Remark \ref{rem:indep}), we see that $K^{S(\mathfrak p)}_*(\Res{G}{S(\mathfrak p)}  A )_{\mathfrak q(\mathfrak p) }$ vanishes if and only if $K^{S}_*(\Res{G}{S}  A )_{\mathfrak r }\simeq 0$ for \emph{some} pair $(S,\mathfrak r)$ with $S$ cyclic and $\mathfrak r\in \spec(R(S))$ lying above $\mathfrak p$. 
By considering all $\mathfrak p$ simultaneously, the above expression becomes
\begin{eqnarray*}
\spec (R(G))\smallsetminus \sigma_G(A) 
&=&
\bigcup_{S} \spec (\Res{G}{S})^{-1}
\big( \spec(R(S))\smallsetminus \Supp_{R(S)}K^S_*(\Res{G}{S}A) \big)
\end{eqnarray*}
where the sum is over all cyclic subgroups of~$G$.
Since $\Res{G}{S}(A)\in \tcat^S_c$, the $R(S)$-module $K^S_*(\Res{G}{S}A)$ is finitely generated (Lemma \ref{lemma:fg}). Therefore its module-theoretic support $\Supp_{R(S)}$ is closed in $\spec R(S)$, and we conclude from the latter formula that $\sigma_G(A)$ is a closed subset of $\spec R(G)$.
\qed

In the next section we prove the last claim of Theorem \ref{thm:gen_supp_intro}.

\subsection{Split injectivity of $f_G: \spec R(G)\to \spc \mathcal K^G$}
\label{subsec:split}

In \cite{balmer_spec2}, Balmer shows that, for every $\otimes$-triangulated category $\tcat$, there is a natural continuous \emph{comparison map}
\begin{equation*}
\rho_{\tcat}:\spc(\tcat)\to \spec(\mathrm R_{\tcat}) \quad, \quad \mathcal P \mapsto \rho_{\tcat}(\mathcal P):=\{r\in \mathrm R_{\tcat} \mid \cone(r)\not\in \mathcal P \}
\end{equation*}
between the spectrum of $\tcat$ and the Zariski spectrum of its central ring. Since the ring $\mathrm R_{\mathcal K^G}=R(G)$ is noetherian (at least for $G$ a compact Lie group), it follows from \cite[Thm.\ 7.3]{balmer_spec2} that $\rho_{\mathcal K^G}: \spc(\mathcal K^G)\to \spec(R(G))$ is surjective. In the previous section, we have constructed a support datum $(\spec(R(G)),\sigma_G)$ on $\mathcal K^G$ for each finite group~$G$. By the universal property of Balmer's spectrum (Prop.\ \ref{prop:UPsupp}), we have the canonical continuous map
\begin{equation*}
f_G:\spec(R(G))\to \spc(\mathcal K^G) \quad , \quad
\mathfrak p \mapsto f_G(\mathfrak p)= \{A\in \mathcal K^G\mid \mathfrak p\not\in \sigma_G(A)\}.
\end{equation*}
We now verify that $f_G$ provides a continuous section of $\rho_{\mathcal K^G}$:

\begin{prop} \label{prop:split}
The composition $\rho_{\mathcal K^G}\circ f_G$ is the identity map of $\spec(R(G))$.
\end{prop}

\proof
Notice that $f_G(\mathfrak p)= \Ker(F_{\mathfrak p})\cap \mathcal K^G$.  
 For a $\mathfrak p\in \spec(R(G))$ and an $r\in R(G)$ we have equivalences (write $\rho:= \rho_{\mathcal K^G}$ and $f:=f_G$ for readability): 
 $r\not\in \rho(f(\mathfrak p))$ 
 $\Leftrightarrow$ $\cone(r) \in f(\mathfrak p)$ (by definition of $\rho$)
 $\Leftrightarrow$ $F_{\mathfrak p}(\cone(r))\simeq 0$ 
  $\Leftrightarrow$ $K^S_*(\Res{G}{S}(\cone(r)))_{\mathfrak q}\simeq 0$, with $\mathfrak q=\mathfrak q(\mathfrak p)$ and $S=S(\mathfrak p)$ (By Lemma \ref{lemma:sigma-cpt}) 
  $\Leftrightarrow$ $K^S_*(\cone(\Res{G}{S}(r)))_{\mathfrak q}\simeq 0$  (because $\Res{G}{S}$ is triangulated)
 $\Leftrightarrow$ $\Res{G}{S}(r)\in (R(S)_{\mathfrak q})^{\times}$.
 
 Thus: $r\not\in \rho(f(\mathfrak p))$ $\Leftrightarrow $ $\Res{G}{S}(r)\in R(S)_{\mathfrak q}^{\times}$. On the other hand, we also have $r\not\in \mathfrak p$ $\Leftrightarrow$ $r\in R(G)_{\mathfrak p}^{\times}$.
 Now observe the commutative square
\begin{equation*}
\xymatrix{
R(G) \ar[d]_{\ell_{\mathfrak p}} \ar[r]^-{\Res{G}{S}} & R(S) \ar[d]^{\ell_{\mathfrak q}} \\
R(G)_{\mathfrak p} \ar[r]^-{} & R(S)_{\mathfrak q}
}
\end{equation*}
where the vertical maps are the localization homomorphism of rings at the indicated prime. Since $\mathfrak p=(\Res{G}{S})^{-1}(\mathfrak q)$, the lower horizontal map is a local homomorphism of local rings, and we deduce that $\ell_{\mathfrak p}(r)$ is invertible if and only if $\ell_{\mathfrak q}(\Res{G}{S}(r))$ is invertible. This proves that $\rho(f(\mathfrak p))=\mathfrak p$.
\qed

\subsection{The spectrum and the Bootstrap category}
\label{subsec:boot}

Theorem \ref{thm:abstract} and Proposition \ref{prop:S7S8} can be easily applied to $\tcat^G=\langle \unit \rangle_{\loc}\subset \KK^G$ in the case of the trivial group, i.e., to the ``Bootstrap category'' $\Boot=\langle \C \rangle_{\loc}\subset \KK$. Its central ring $R(G)$ is just $\Z$, and its subcategory of compact objects $\Boot_c=\langle \C\rangle$ is the full subcategory of separable $\Cstar$-algebras having finitely generated $K$-theory groups (see \cite[Lemma 5.1.6]{thesis}).

\begin{thm} \label{thm:spec_bootc}
There is a canonical isomorphism $\spec(\Boot_c)\simeq \spec(\Z)$ of locally ringed spaces, given by $\rho_{\Boot_c}$ with inverse $f_G$.
\end{thm}

\proof
Let $\sigma: \obj(\Boot)\to 2^{\spec (\Z)}$ be the support constructed in $\S$\ref{subsec:supp_datum}, for $G=\{1\}$. 
Namely: $\sigma(A)= \{(p)\in \spec(\Z)\mid \mathbb F_p\otimes_{\Z}^{\mathsf L} K_*(A) \not\simeq 0\}$ (here $\mathbb F_0:=\Q$). In this case at least, $\sigma$ detects objects (see \cite[Lemma 2.12]{neemanChr} for a more general statement working for any commutative noetherian ring $R$ instead of~$\Z$). Moreover, if $A\in \Boot_c$ then $\sigma(A)=\{(p)\mid K_*(A)_{(p)} \not\simeq 0\}=\Supp_{\Z}(K_*(A))$   by Lemma \ref{lemma:sigma-cpt}.
Thus, by Theorem \ref{thm:sigma} and Proposition \ref{prop:S7S8}, $\sigma$ satisfies \emph{all} ten hypotheses (S0)-(S9) of Theorem \ref{thm:abstract}, and therefore we have a canonical homeomorphism $f:=f_{\{1\}}:\spc(\Boot_c)\simeq \spec(\Z)$.
By Proposition \ref{prop:split}, its inverse must be the comparison map $\rho:=\rho_{\Boot_c}$.
It is now a general fact, true for any $\otimes$-triangulated category $\tcat$, that if $\rho_{\tcat}$ is a homeomorphism then it yields also automatically an isomorphism of locally ringed spaces $\spec(\tcat)\simeq \spec(\mathrm R_{\tcat})$; see \cite[Prop.\ 6.11 (b)]{balmer_spec2}. Alternatively, in the case at hand it is straightforward to check this directly.
\qed

\begin{remark}
In \cite[$\S$5.1]{thesis} we give a more elementary proof of Theorem \ref{thm:spec_bootc}, relying on the classical Universal Coefficient theorem and the K\"unneth theorem of Rosenberg and Schochet \cite{rs}.
\end{remark}


\end{document}